\documentclass{amsart}
\usepackage{amssymb,amsmath,amscd,amsfonts}
\usepackage[all]{xy}
\newcommand{\commentout}[1]{}

\newtheorem{Theorem}{Theorem}[subsection]
\newtheorem{Proposition}[Theorem]{Proposition}
\newtheorem{Corollary}[Theorem]{Corollary}
\newtheorem{Lemma}[Theorem]{Lemma}

\newtheorem{Remark}[Theorem]{Remark}
\newtheorem{Definition}[Theorem]{Definition} 
\newtheorem{Exercise}[Theorem]{\exercisename}
\newtheorem{Example}[Theorem]{Example}

\newtheorem{Conjecture}[Theorem]{Conjecture}
\newtheorem{Problem}[Theorem]{Problem}
\newtheorem{Problems}[Theorem]{\problemsname}

\newenvironment{cor}{\begin{Corollary}}{\end{Corollary}}
\newenvironment{lem}{\begin{Lemma}}{\end{Lemma}}
\newenvironment{re}{\begin{Remark}}{\end{Remark}}
\newenvironment{defn}{\begin{Definition}}{\end{Definition}}
\newenvironment{ex}{\begin{Example}}{\end{Example}}

\newcommand{\mc}[1]{\mathcal#1}
\newcommand{\h}[1]{{\widehat{#1}}}
\newcommand{\f}[1]{{\mathfrak{#1}}}
\newcommand{\diag}{{\rm diag}}

\newcommand{\Hom}{{\rm Hom}}
\newcommand{\ol}{\overline}

\newcommand{\mb}[1]{\mathbb#1}
\newcommand{\Ind}{{\rm Ind}}
\begin{document}
\pagenumbering{arabic}
\title{Certain Unitary Langlands-Vogan Parameters for Special Orthogonal Groups}
\author{Hongyu HE} 
\footnote{Email address: hongyu@math.lsu.edu; Fax: 225-578-4276.}
\footnote{This research is partially supported by the NSF grant DMS 0700809.}
\address{Department of Mathematics, Louisiana State University,
Baton Rouge, LA 70803}
\begin{abstract}
\vspace{0.5 in}
{\large
A Langlands parameter, in the Langlands dual group, can be decomposed into a product of a tempered parameter and a positive quasi-character. Fixing a tempered parameter, Arthur  conjectured that positive quasi-characters corresponding to certain weighted Dynkin diagrams for the centralizer of the tempered parameter will yield unitary representations. In this paper, we treat one basic case in Arthur's conjecture for the special orthogonal groups. We establish the unitarity for a class of Langlands-Vogan parameters in Arthur's packets. }
\vspace{0.5 in}
 \end{abstract}

\maketitle
\clearpage
\section*{Introduction}
Let $G$ be a semisimple Lie group. Let $\Pi_u(G)$ be the unitary dual of $G$. Let $\Pi(G)$ be the admissible dual of $G$. Fix a minimal parabolic subgroup $P_m$. Let $P$ be a parabolic subgroup  containing $P_m$. Let $MAN$ be the Langlands decomposition of $P$. Langlands showed that $\Pi(G)$ is in one-to-one correspondence with a triple $(P, \sigma, v)$, where $\sigma$ is the (infinitesimal) equivalence class of an irreducible tempered representation of $M$ and $v$ is a complex-valued linear functional on $\f a$ such that its real part $\Re(v)$ is in the open positive Weyl chamber of $\f a^*$. The representation corresponding to $(P, \sigma, v)$ is $J(P, \sigma, v)$, the Langlands quotient. By a theorem of Harish-Chandra, $\Pi_u(G)$ can be regarded as a subset of $\Pi(G)$. 
We shall say that a Langlands parameter $(P, \sigma, v)$ is unitary or unitarizable if $J(P, \sigma, v)$ is unitarizable. \\
\\
Irreducible tempered representations are all unitarizable. Therefore their equivalence classes, the tempered dual, can be regarded as a subset of $\Pi_u(G)$. The classification of the tempered dual is completed by Knapp and Zuckerman (\cite{kz}). Fix a tempered parameter $\sigma$, the problem of determining the unitary dual $\Pi_u(G)$ can be approached by determining those $v$ such that $J(P, \sigma, v)$ is unitarizable.  This set of $v$ is necessarily bounded.  In \cite{ar}, Arthur conjectured that if $v$ corresponds to certain weighted Dynkin diagram, $J(P, \sigma, v)$ will be unitary. In this paper, we shall treat one basic case in Arthur's conjecture for special orthogonal groups.  \\
\\
Before we state our main results. Let me introduce some notations. Let $G=SO(p,q)$. Suppose that $0 < p \leq q$. Fix a maximal compact subgroup $K$. Then $K \cong S(O(p) O(q))$. $K$ is disconnected.  We parametrize the unitary dual of $SO$ by dominant weights. Let  $(\xi) \in \Pi_u(SO(p))$ and $(\eta) \in \Pi_u(SO(q))$ as in \ref{compactrepn}. If $(\xi) \otimes (\eta)$ extends to a representation of $K$, then we obtain two such extensions $(\xi, \eta, \pm) \in \Pi_u(K)$. We use the convention set in \cite{kv} to mark these two different extensions. If $(\xi) \otimes (\eta)$ does not extend to a representation of $K$, then there is a unique representation $(|\xi|, |\eta| ) \in \Pi_u(K)$ such that $(|\xi|, |\eta|)|_{SO(p)SO(q)}$ contains a copy of $(\xi) \otimes (\eta)$. In order to make our result easier to state, we denote the representation $(|\xi|, |\eta|)$ by $(|\xi|, |\eta|, +)$. The unitary dual of $S(O(p) O(q))$ can be parametrized by certain $(\xi, \eta, \pm)$ and $(|\xi|, |\eta|, +)$. If $p=0$, we will write $(\eta)$ as $(0, \eta, +)$. See \ref{compactrepn} for details. \\
\\
Given a {\bf real} vector $\mu=(\mu_1, \mu_2, \ldots, \mu_n)$, let
$\overline{\mu}$ be the rearrangement of $\mu$ in descending order
$$\overline{\mu_1} \geq \overline{\mu_2} \geq \ldots \geq \overline{\mu_n}.$$
Let $| \mu |=(| \mu_1 |, | \mu_2 |, \ldots | \mu_n |)$. For any constant $C$, let $\mathbf c$ be the constant vector $(c, c, \ldots, c)$ of a proper dimension. All general linear groups in this paper are real general linear groups, denoted by $GL(n)$. 
\begin{Theorem} Let $\sigma$ be an irreducible tempered representation of $SO(p-d, q-d)$ such that one of its minimal $S(O(p-d)O(q-d))$-types is $(\xi, \eta, +)$ (\cite{vogan79}). Then the Langlands-Vogan subquotient of 
$\Ind_{SO(p-d, q-d) GL(1)^d N}^{SO(p,q)}  \sigma \otimes \prod_{i=1}^d |\det|^{v_i}$ is unitary where
\begin{enumerate}
\item $v$ is the semisimple element in $\f {so}(2d, \mathbb C)$ corresponding to $\diag(\frac{1}{2}, -\frac{1}{2})$ under a Lie algebra homomorphism  $\f {sl}(2, \mathbb C) \rightarrow \f {so}(2d, \mathbb C)$ when $p+q$ is even;
\item $v$ is the semisimple element in $\f {sp}(2d, \mathbb C)$ corresponding to $\diag(\frac{1}{2}, -\frac{1}{2})$ under a Lie algebra homomorphism $\f {sl}(2, \mathbb C) \rightarrow \f {sp}(2d, \mathbb C)$ when $p+q$ is odd.
\end{enumerate}
Here $v$ is parametrized by $\mathbb R^d$ and $|\det|$ is simply the absolute value of $g \in GL(1)$.
\end{Theorem} 
If $d \neq 0$, the Langlands-Vogan subquotient is defined to be the irreducible subquotient  of 
$$\Ind_{SO(p-d, q-d) GL(1)^d N}^{SO(p,q)}  \sigma \otimes \prod_{i=1}^d |\det|^{v_i}$$ containing a minimal $K$-type $(|\xi| \oplus {\mathbf 0}, |\eta| \oplus {\mathbf 0}, +)$. If $(p,q) \neq (d,d)$, it is equivalent to the Langlands quotient of 
$$\Ind_{SO(p-d, q-d) GL(1)^d N}^{SO(p,q)}  \sigma \otimes \prod_{i=1}^d |\det|^{|\overline{v_i}|}.$$ When $\sigma$ is trivial,  we obtain unipotent representations. In this case, our result overlaps with \cite{baru} (for $p=q, q-1$) and \cite{quan}. \\
\\
Originally, Arthur's conjecture was  formulated for quasisplit groups (\cite{ar}). For quasisplit groups, each Arthur's packet yields a Langlands' packet. The passage from Arthur's packet to Langlands' packet is quite natural and elegant. Arthur's conjecture can thus be formulated directly in terms of Langlands parameters. For nonquasisplit groups, the formulation of Arthur's packets is quite complicated and the precise formulation of Arthur's conjecture can be found in \cite{abv}. However, even for nonquasisplit groups, Arthur's conjecture, in its original form (\cite{ar}), still makes sense and can still be formulated in terms of Langlands parameter. Our main theorem confirms a basic case of Arthur's conjecture under this
rather restrictive interpretation.   
The basic idea is to treat $J(\sigma, v)$ uniformly for all tempered parameter $\sigma$. In particular, we identify $J(\sigma, v)$  with a certain quantum induced module, which appears as an invariant tensor product
(Definitions \ref{average} \ref{itp} \ref{quantuminduction}). The unitarity of this quantum induced module can then be established by some standard argument. \\
\\
This paper is organized as follows. In Section 1, we define the invariant tensor product of two representations weakly as distributions (Definitions \ref{average} \ref{itp}). We prove that invariant tensor product may inherit Lie group actions (Lemma \ref{gkh}), Lie algebra actions (Lemma \ref{tensorrepn1}) and Hermitian forms (Prop. \ref{invher}). In Section 2, we review  the basic theory of induced representations, Langlands classification, Vogan's subquotients and  growth of the matrix coefficients.
In Section 3, we review the known results on the degenerate principal series $I_n(s)$ and its unitarizable subquotients, due to Johnson (\cite{jo}), later reworked by Sahi (\cite{sa}), Zhang (\cite{zh}) and Lee (\cite{ll}). We give some additional analysis on the small constituents $\mc E_m(n)$ with $0 \leq m \leq [\frac{n-1}{2}]$. \\
\\
In Section 4, we define the quantum induction $Q(2m)(\pi)$ on the Harish-Chandra module level to be $V(\mc E_m(p+q+d)) \otimes_{SO(p,q)} V(\pi)$. Here $V(\cdot)$ stands for the Harish-Chandra module and $\pi \in \Pi(SO(p,q))$ satisfies a certain growth condition. $Q(2m)(\pi)$ is a Harish-Chandra module of $SO(q+d, p+d)$.
When $\pi$ is unitary, $\mc Q(2m)(\pi)$ inherits an invariant Hermitian structure from $\mc E_m(p+q+d)$. Suppose that $p+q \leq 2m+1 \leq p+q+d$.  The main task in Section 4 is to show that the canonical invariant Hermitian form for $\mc Q(2m)(\pi)$ is positive definite under a growth condition on $\pi$ (Theorem \ref{posinv}). \\
\\
In Section 5, we show that $\mc Q(2m)(\pi)$ is a subrepresentation of
$\Ind_{SO(q,p) GL(d) N}^{SO(q+d, p+d)} \pi \otimes |\det|^{m-\frac{n-1}{2}}$ if we regard $\pi$ as a representation of $SO(q,p)$ (Theorem \ref{quan_ind}). \\
\\
In Section 6, we relate $\mc Q(2m)$ to the composition of Howe's correspondence (\cite{howe}). Some caution is taken to handle the difference between $SO(p,q)$ and $O(p,q)$. Suppose that $p+q \leq 2m+1 \leq p+q+d$. We apply Kudla's preservation principle and results of Moeglin and Adam-Barbasch to show the nonvanishing of $\mc Q(2m)(\pi)$ (\cite{ku} \cite{mo} \cite{ab} \cite{non}).  See Theorem \ref{quan_theta}.  We then apply Howe's theory (\cite{howe} and results of \cite{mo} \cite{ab} \cite{paul} to determine a minimal $K$-type of $\mc Q(2m)(\pi)$. Our assumption that $\pi$ has a $K$-type of the form $(\xi, \eta, +)$ is essential to guarantee that the minimal $K$-type is a $K$-type of minimal degree. This allows us to determine the Langlands-Vogan parameter of $\mc Q(2m)(\pi)$. See Theorem \ref{quan_vogan}. \\
\\
Finally, in Section 7, we review some basic facts about Arthur's packets and show that quantum induction can be applied inductively to obtain unitary Langlands-Vogan parameters in Arthur's packets. One highlight is the following theorem.
\begin{Theorem}\label{0.2} Suppose that $p+q \leq 2m+1 \leq p+q+d$ and $p \leq q$. Let $\pi$ be an irreducible unitary representation of $SO(q,p)$ such that
its every $K$ finite matrix coefficient $f(g)$ satisfies the condition that
$$ |f(k_1 \exp H(g) k_2)| \leq C_f \exp (\overbrace{m+2-p-q- \epsilon, m+3-p-q, \ldots m+1-q}^{p})(| \overline{H(g)}|)$$
for some $\epsilon >0$. Here $k_1, k_2 \in S(O(q)O(p))$ and $H(g) \in \f a_{\f m} \cong \mathbb R^p$. Then there exists a unitaizable subrepresentation of
$\Ind_{SO(q,p) GL(d) N}^{SO(q+d, p+d)} \pi \otimes \mathbb |\det|^{m-\frac{p+q+d-1}{2}}$.
\end{Theorem}
Here $p \leq q$ is not essential. See Theorem A in Section 7 for more details. I shall remark that Theorem \ref{0.2} can be applied to construct many more unitary representations. However, at this point, we are unable to identify the Langlands-Vogan parameters of these unitary representations.

\subsection{Notations}
Write $ \mu \prec \lambda$ if for every $k \in [1, n]$
$$\sum_{i=1}^{k} \mu_i < \sum_{i_1}^{k} \lambda_i.$$
Write $\mu \preceq \lambda$ if for every $k \in [1, n]$
$$\sum_{i=1}^{k} \mu_i \leq \sum_{i_1}^{k} \lambda_i.$$
Let $\mathbb Z$ be the set of integers. Let $\mathbb Z +\frac{1}{2}$ be the set of half integers.\\
\\
 All topological groups and topological vector spaces are assumed to be Hausdorff. Let $\mc X$ be a topological vector space (TVS). Let $\mathcal L(\mathcal X)$ be the space of continuous linear operators on  $\mc X$. A linear representation $\pi$  of $G$ on $\mathcal X$, is said to be continuous 
if the group action
$ G \times \mathcal X \rightarrow \mathcal X$ is
continuous.  The TVS $\mc X$ here needs not to be complete. \\
\\
If $\mc X$ is a Hilbert space, then $(\pi, \mc X)$ is called a Hilbert representation. 
All Hilbert spaces in this paper are assumed to be separable. All Hilbert representations are assumed to be continuous. \\
\\
For any complex vector space $V$, let $V^c$ be the vector space $V$, regarded as a real vector space, equipped with the conjugated complex  multiplication, 
$$ c. u \in V^c = \overline{c} u \in V \qquad (c \in \mathbb C, u \in V^c).$$
Then $V^c$ is a complex linear vector space. $V^c=V$ as real vector spaces. But $V^c$ and $V$ have different complex structures.\\
\\
A maximal compact subgroup of a semisimple Lie group may be universally denoted by $K$.   The nilradical of a parabolic subgroup of a semisimple group will be universally denoted by $N$. It should be clear within the context what $K$ and $N$ are. $C$ will be used as a universal constant. Identity operators or matrices will be denoted by $I$ or $I_n$ where $n$ is the dimension. \\
\\
All the important conventions will be highlighted by boldfaced letters. All inner products or Hermitian forms will be denoted by $(\, , \,)$. The definition of these forms will be clear within the context.\\
\\
Let $G$ be a unimodular Lie group. Let $X$ be a principal $G$-bundle with a $G$-invariant measure. Let $\mc H$ be a unitary representation of $G$. Then $L^2(X \times_G \mc H, X/G)$ is a Hilbert space. The Hilbert inner product is given by
$$(\phi, \psi)=\int_{[x] \in X/G} (\phi(x), \psi(x)) d [x].$$
Whenever we have an integral of this form, the integrant does not depend on the choices of $x$ in $[x]$.
\section{Invariant Distributions and Invariant Tensor Product} 
\subsection{Invariant Distributions and Averaging Operators}
Let $G$ be a locally compact  group. Let $(\pi, X)$ be a linear representation of $G$. Let $u \in X$. $u$ is called a $G$-invariant vector if $\pi(g)(u)=u$ for any $g \in G$. We denote the space of invariants of $ X$ by
$ X^G$.\\
\\
Let $\pi$ be a continuous representation of $G$ on a locally convex TVS $\mc X$. 
Let $\mc X^*$ be the dual space of $\mc X$ equipped with the weak-* topology. The continuous representation $(\pi, \mathcal X)$ induces a linear representation $(\pi^{*}, \mathcal X^{*})$ as follows. For any $\delta \in \mathcal X^{*}$, $v \in \mathcal X$, $g \in G$, define
$$\pi^{*}(g)(\delta)(v)=\delta(\pi(g)^{-1} v).$$
This is the dual representation.
Obviously, $\pi^{*}(g)(\delta) \in \mathcal X^{*}$. However, the dual representation
$G \times \mathcal X^{*} \rightarrow \mathcal X^{*}$
may not be continuous.\\
\\
Given $\delta \in \mathcal X^{*}$ and $v \in \mathcal X$, the {\bf matrix coefficient}
$$ g \rightarrow \delta(\pi(g) v)$$
is a continuous function on $G$. Denote it by $\pi_{v, \delta}(g)$.  
\begin{defn}\label{average}
Let $G$ be a locally compact unimodular group with finite center.
Suppose that there is a subspace $\mathcal Y $ of $\mathcal X^{*}$, and a subspace $\mathcal X_0$ of $\mathcal X$ such that $\pi_{v, \delta}(g) \in L^1(G)$ for any $\delta \in \mathcal Y$ and
any $v \in \mathcal X_0$. Then we define a map
$$\mathcal A_{G, \mathcal Y}: \mathcal X_0 \rightarrow  \Hom(\mathcal Y, \mathbb C)$$
by $\, \forall \,\, u \in \mathcal X_0$,
$$\mathcal A_{G, \mathcal Y}(u)(\delta)=\int_{G} \delta(\pi(g) u) d g, \qquad (\delta \in \mathcal Y).$$
We call $\mathcal A_{G, \mathcal Y}$ the averaging operator with respect to $\mathcal Y$.
\end{defn}
\begin{lem}
 If $\mathcal Y$ is $G$-stable, then the image of the averaging operator
$$\mathcal A_{G, \mathcal Y}(\mathcal X_0) \subseteq \Hom(\mathcal Y, \mathbb C)^G.$$
\end{lem}
Proof: For any $u \in \mc X_0, \delta \in \mc Y, g \in G$, we have
\begin{equation}
\begin{split}
 A_{G, \mc Y}(u)(\pi^*(g) \delta) = & \int_{h \in G} (\pi^*(g) \delta)(\pi(h)u) d h \\
 = &\int_{h \in G}\delta(\pi(g^{-1}) \pi(h)u) d h \\
 = & \int_{h \in G} \delta(\pi(h)u) d h \\
 = & \mc A_{G, \mc Y}(u)(\delta).
 \end{split}
 \end{equation}
We shall remark that whether $\mc Y$ is $G$-stable or not is not essential. Consider 
$${}^G\mc Y={\rm span}\{ \pi^{*}(g) (\delta) \mid g \in G, \delta \in \mc Y \}.$$
Then $\mc A_{G, {}^G\mc Y}(\mc X_0)$ is well-defined and 
$$\mc A_{G, {}^G\mc Y}(\mc X_0) \subseteq \Hom({}^G\mathcal Y, \mathbb C)^G.$$

\subsection{Choices of $\mc Y$ and $L^1$-mutually Dominated Subspaces}
Natually, $\mathcal A_{G, \mathcal Y}(\mathcal X_0)$ shall be interpreted as a linear subspace of linear functions on $\mathcal Y$. The averaging operator will then depend on the choices of the subspace $\mathcal Y$ in $\mathcal X^{*}$. This is inconvenient in applications.  There is  an interpretation of $\mathcal A_{G}$ that is somewhat independent of the choices of $\mathcal Y$ as I shall explain.\\
\\
Notice that $$\mathcal A_{G, \mathcal Y}(\mathcal X_0) \cong \mathcal X_0/ker(\mathcal A_{G, \mathcal Y}).$$
If for two different choices $\mathcal Y_1$ and $\mathcal Y_2$, we can show that 
$$ker(\mathcal A_{G, \mathcal Y_1})=ker(\mathcal A_{G, \mathcal Y_2}),$$
then the image $\mathcal A_{G, \mathcal Y_1}(\mathcal X_0)$ can be identified with $\mathcal A_{G, \mathcal Y_2}(\mathcal X_0)$.

\begin{Definition}\label{l1dominated} Fix $\mathcal X_0 \subseteq \mathcal X$. Let $\mc Y_1, \mc Y_2$ be two  subspaces of $\mathcal X^{*}$. We say that  $\mathcal Y_1$ is $L^1$-dominated by $ \mathcal Y_2$ (with respect to $\mathcal X_0$),  if for every $u \in \mathcal X_0$ and $\delta \in \mathcal Y_1$, there is a sequence $\{ f_{\alpha} \} \subset \mathcal Y_2$ such that
$$\pi_{u, f_{\alpha}}(g) \rightarrow \pi_{u, \delta}(g) \qquad (\forall \, \, g \in G)$$  and $\{ \pi_{u, f_{\alpha}}(g) \}$ are uniformly dominated by an $L^1$-function on $G$. If $\mathcal Y_2$ is also $L^1$-dominated by $\mathcal Y_1$, we say that $\mc Y_1$ and $\mc Y_2$ are mutually $L^1$-dominated subspaces of $\mathcal X^{*}$ with respect to $\mc X_0$.
\end{Definition}

\begin{Remark}\label{l1dominatedremark}
If $\mc Y_1$ are the (finite) linear span of a set of vectors, it suffices that
the condition in Definition \ref{l1dominated} holds for $\delta$ in the spanning set.
Same is true for a spanning set of vectors in $\mc X_0$.
\end{Remark}

\begin{Theorem}\label{l1dominatedTheorem} Suppose that $\mathcal Y_1$ and $\mathcal Y_2$ are mutually $L^1$-dominated subspaces of $\mathcal X^{*}$ with respect to $\mc X_0$. Then
$$ker(\mathcal A_{G, \mathcal Y_1})=ker(\mathcal A_{G, \mathcal Y_2}),$$
and
$\mathcal A_{G, \mathcal Y_i}(\mathcal X_0)$ can be identified with each other cannonically. 
\end{Theorem}
Proof: Suppose $u \in ker(\mathcal A_{G, \mathcal Y_2}) \subset \mathcal X_0$. Then for any $f \in \mathcal Y_2$, we have
$$\mathcal A_{G, \mathcal Y}(u)(f)=\int_{G} f(\pi(g) u) d g=\int \pi_{u, f}(g) d g=0.$$
Let $\delta$ be arbitrary in $ \mathcal Y_1$.
Since $\mathcal Y_1$ is $L^1$-dominated by $ \mathcal Y_2$, there is a sequence $\{ f_{\alpha} \} \subset \mathcal Y_2$ such that
$\pi_{u, f_{\alpha}}(g) \rightarrow \pi_{u, \delta}(g)$ and $\pi_{u, f_{\alpha}}(g)$ are uniformly dominated by an $L^1$-function on $G$.
Since $\int_G \pi_{u, f_{\alpha}}(g) d g=0$, by the dominated convergence theorem, $\int_G \pi_{u, \delta}(g) d g=0$. Hence $ u \in ker(\mathcal A_{G, \mathcal Y_1})$. We have shown that
$$ker(\mathcal A_{G, \mathcal Y_2}) \subseteq ker(\mathcal A_{G, \mathcal Y_1}).$$
The converse can be proved in the same way. Our assertions follow. $\Box$  
\subsection{Averaging Operator $\mathcal A_G$ for Hermitian Representations }
\begin{defn}
Let $(\pi, \mathcal X)$ be a continuous representation of $G$ on a TVS $\mc X$.  If $\mc X$ is equipped with a nondegenerate Hermitian form, we say that $(\pi, \mc X)$ is a {\it Hermitian representation} of $G$. Let $V$ be a subspace of $\mc X$ such that for any $v_1, v_2 \in V$,
$$\pi_{v_1, v_2}: g \in G \rightarrow (\pi(g) v_1, v_2)$$
is in $L^1(G)$. Define
$$\mathcal A_G: V \rightarrow \Hom(V^c, \mathbb C)$$
by
$$\mathcal A_G(v_1)(v_2)=\int_G (\pi(g)v_1, v_2) d g \qquad (v_1 \in V, v_2 \in V^c).$$
\end{defn}

Let $(\pi, \mathcal H)$ be a (continuous) Hilbert representation of $G$.  Let $(\pi^{*}, \mathcal H^{*})$ be the contragredient Hilbert representation.
Let $\mathcal H^c$ be $\mathcal H$, equipped with the conjugate complex structure
$$ c. u \in \mathcal H^c= \overline{c} u \in \mathcal H \qquad (c \in \mathbb C).$$
Clearly, $\pi(g)$ acts on $\mathcal H^c$ preserving the complex structure.
So it defines a Hilbert  representation, denoted by $(\pi^c, \mathcal H^c)$. If $\pi$ is unitary, we have
$$ (\pi^c, \mathcal H^c) \cong (\pi^{*}, \mathcal H^{*}),$$
essentially by the Riesz representation theorem.
\begin{re} Let $(\pi, \mathcal H)$ be a (continuous) Hilbert representation of $G$.
Suppose that there is a subspace $V \subseteq \mathcal H$ such that
for any $v_1, v_2 \in V$,
$$\pi_{v_1, v_2}: g \in G \rightarrow (\pi(g) v_1, v_2)$$
is in $L^1(G)$. Then 
$\mathcal A_G: V \rightarrow \Hom(V^c, \mathbb C)$ is well-defined. On the other hand, by Riesz representation theorem, the inner product on $\mc H$ induces a topological embedding
$ V^c \hookrightarrow \mc H^*$.  The canonical action of $G$ on $V^c$  is $\pi^*$.
Thus $\mathcal A_{G, V^c}$ is well-defined and  $\mc A_{G}(V)=\mc A_{G, V^c}(V)$. 
{\bf Unless stated otherwise, $V^c$ will be equipped with the action of $\pi^*$}.
\end{re}
If $V$ is $G$-invariant, then $\mc A_{G}(V)=\mc A_{G, V^c}(V)$ will be a subspace of $\Hom(V^c, \mathbb C)$.
Generally speaking $V^c$ will not be $\pi^*(G)$-invariant unless
$(\pi|_V, V)$ is unitary. So the image $\mc A_{G}(V)$ may not be in the $G$-invariant subspace of $\Hom(V^c, \mathbb C)$.
\begin{ex}
Let $G$ be a compact group equipped with the probability measure. Let $(\pi, \mathcal H)$ be a unitary representation of $G$. Then $\mathcal A_G$ is defined for the whole space $\mathcal H$. The space 
$$\mathcal A_G(\mathcal H) = \Hom(\mathcal H^c, \mathbb C)^G.$$
The right hand side can be identified with $\mathcal H^G$. The operator $\mathcal A_G$ is the projection operator onto the trivial-isotypic subspace. Consequently, for $G$ compact, the averaging integral
 $\mc A_G(v)$ can simply be defined as $\int_G \pi(g) v d g$. This can be regarded as the strong form of the averaging operator.

\end{ex}

\subsection{Invariant Tensor Products}
\begin{defn}\label{itp}
Let $G$ be a unimodular group. Let $(\pi, \mathcal X)$ and $(\tau, \mathcal Z)$ be two continuous representations of $G$. Suppose that there are  subspaces
$\mathcal X_0 \subseteq \mathcal X,  \mathcal Z_0 \subseteq \mathcal Z$ and  subspaces $ \mathcal Y \subseteq \mathcal X^{*},  \mathcal W \subseteq \mathcal Z^{*}$
such that the matrix coefficents
$$\pi_{\mathcal X_0, \mathcal Y}(g) \tau_{\mathcal Z_0, \mathcal W}(g) \subseteq L^1(G).$$
We define
$$ \mathcal X_0 \otimes_{G; \mathcal Y \otimes \mathcal W} \mathcal Z_0=\mathcal A_{G, \mathcal Y \otimes \mathcal W}(\mc X_0 \otimes \mc Z_0) \subseteq \Hom(\mathcal Y \otimes \mathcal W, \mathbb C).$$
We call $\mathcal X_0 \otimes_{G; \mathcal Y \otimes \mathcal W} \mathcal Z_0$ the invariant tensor product of $\mathcal X_0$ and $\mathcal Z_0$ with respect to $\mathcal Y \otimes \mathcal W$. 
\end{defn}
Similarly, we can define $\mc Y \otimes_{G, \mc X_0 \otimes \mc Z_0} \mc W$ since $\mc X_0 \subseteq \mc X \subseteq (\mc X^*)^*$ and $\mc Z_0 \subseteq \mc Z \subseteq (\mc Z^*)^*$.
\begin{lem}\label{pairing}
For any $x_i \in \mc X_0, z_i \in \mc Z_0, y_j \in \mc Y, w_j \in \mc W$,
$$(\sum_{i=1}^t x_i \otimes_{G, \mc Y \otimes \mc W} z_i)(\sum_{j=1}^s y_j \otimes w_j)=
(\sum_{j=1}^s y_j \otimes_{G, \mc X_0 \otimes \mc Z_0} w_j)(\sum_{i=1}^t x_i \otimes z_i).$$
This induces a nondegenerate bilinear form
$$ (\mc X_0 \otimes_{G, \mc Y \otimes \mc W} \mc Z_0, \mc Y \otimes_{G, \mc X_0 \otimes \mc Z_0} \mc W) \rightarrow \mathbb C,$$
namely 
$$(\sum_{i=1}^t x_i \otimes_{G, \mc Y \otimes \mc W} z_i, \sum_{j=1}^s y_j \otimes_{G, \mc X_0 \otimes \mc Z_0} w_j)=(\sum_{i=1}^t x_i \otimes_{G, \mc Y \otimes \mc W} z_i)(\sum_{j=1}^s y_j \otimes w_j).$$
\end{lem}
Proof: The first equation follows directly from definition:
$$(\sum_{i=1}^t x_i \otimes_{G, \mc Y \otimes \mc W} z_i)(\sum_{j=1}^s y_j \otimes w_j)=\int_{G} \sum_{j=1}^s \sum_{i=1}^t
y_j (\pi(g) x_i) w_j(\tau(g) z_i) d g.$$
Notice that
$(\sum_{i=1}^t x_i \otimes_{G, \mc Y \otimes \mc W} z_i)(\sum_{j=1}^s y_j \otimes w_j)=0$ for all $\{ x_i \}
\subseteq \mc X_0, \{z_i \} \subseteq \mc Z_0$ if and only if
$\mc A_{G, \mc X_0 \otimes \mc Z_0}(\sum_{j=1}^s y_j \otimes w_j)=0$. So one obtains a form
$$(\sum_{i=1}^t x_i \otimes_{G, \mc Y \otimes \mc W} z_i, \sum_{j=1}^s y_j \otimes_{G, \mc X_0 \otimes \mc Z_0} w_j)=(\sum_{i=1}^t x_i \otimes_{G, \mc Y \otimes \mc W} z_i)(\sum_{j=1}^s y_j \otimes w_j).$$
It is easy to see that this form is nondegenerate. $\Box$.

\begin{defn}
If $(\pi, \mc X)$ and $(\tau, \mc Z)$ are Hermitian representations
and $\mc X_0 \subseteq \mc X$, $\mc Z_0 \subseteq \mc Z$ such that
$$(\pi(g) x, y) (\tau(g) z, w) \in L^1(G) \qquad (x,y \in \mc X_0, z,w \in \mc Z_0),$$ we define the cannonical
$$ \mc X_0 \otimes_G \mc Z_0=\mc A_G(\mc X_0 \otimes \mc Z_0) \subseteq \Hom(\mc X_0^c \otimes \mc Z_0^c, \mathbb C),$$
and
$$ (x \otimes_G z)(y \otimes w)=\int_G (\pi(g)x, y)(\tau(g) z, w) d g, \qquad
(x, y \in \mc X_0; y, w \in \mc Z_0).$$
\end{defn}
{\bf Whenever we use the notation $ \mc X_0 \otimes_G \mc Z_0$, we assume that $G$ is unimodular and $ \mc X_0 \otimes_G \mc Z_0$ is well-defined.}\\
\\
In the case that $(\pi, \mc X)$ and $(\tau, \mc Z)$ are Hilbert representations, then $\mc X_0^c \hookrightarrow \mc X^*$ and $\mc Z_0^c \hookrightarrow \mc Z^*$. We have
$$\mc X_0 \otimes_G \mc Z_0=\mc X_0 \otimes_{G, \mc X_0^c \otimes \mc Z_0^c} \mc Z_0.$$
\commentout{Our definition of invariant tensor product can be modified to include the examples from integration operators. For example, let $f$ be a measuable function on the product of two measured space $X$ and $ Y$, and $p$  a measurable function on $Y$. Suppose that
$f(x,y)p(y)$ is in $L^1(Y)$ for almost every $x$. Then we can define
$$(f \otimes_Y p)(x)=\int_{Y} f(x, y) p(y) d y.$$}

\subsection{A Strong form of ITP---the Geometric Realization}
\begin{Theorem}\label{georeal}
Let $G$ be a unimodular Lie group and $X$ be a smooth manifold. Let $G$ act on $X$ smoothly and freely from the right.  Suppose that $X$ has a nontrivial $G$-invariant  measure $d x$ given by a smooth nowhere vanishing  volume form. Denote the action of $G$ on $L^2(X, dx)$ from the right by $R$. Then $L^2(X)$ is a unitary representation of $G$. Let $C_c^{\infty}(X)$ be the space of smooth function on $X$ with compact support.  Let $(\pi, \mathcal H)$ be a unitary representation of $G$. Then there is an injection $\mc I$ from 
$C_c^{\infty}(X) \otimes_G \mathcal H$ to $ C_c^{\infty}(X \times_G \mc H)$ ( the smooth and compactly supported sections of the Hilbert vector bundle
$X \times_{G} \mathcal H \rightarrow X/G$)  defined by
$$\mc I(f \otimes u)(x)=\int_G f(x g) \pi(g)u d g, \qquad (f \in C_c^{\infty}(X), u \in \mc H).$$
In addition, for any dense subspace $V$ of $\mc H$, $\mc I(C_c^{\infty}(X) \otimes_{G} V)$ is dense in $L^2(X \times_{G} \mc H)$.
\end{Theorem}
\begin{re} Regard $X$ as a principal $G$-fiberation. This theorem says that $\mc A_{G}$ can be interpreted as integration along the fiber in a proper sense. This is a strong form of invariant tensor product.
\end{re}
Proof: Let $f_1, f_2 \in C_c^{\infty}(X)$ and $u_1, u_2 \in \mathcal H$. Then
$(R(g)f_1, f_2)_{L^2(X)}$ must also have compact support on $G$. We have
\begin{equation}\label{gr}
\begin{split}
(f_1 \otimes_G u_1)(f_2 \otimes u_2) =& \int_{G} [\int_X f_1(x g) \overline{f_2}(x) d x] (\pi(g) u_1, u_2) d g\\
= & \int_{G \times X} f_1(x g) \overline{f_2}(x)  (\pi(g) u_1, u_2) d x d g \\
=& \int_{g \in G} \int_{X/G} \int_{h \in G} f_1(x h g) \overline{f_2}(x h) d h d [x] (\pi(g) u_1, u_2) d g \\
=& \int_{X/G} \int_{G \times G} f_1(x h g) \overline{f_2}(x h) (\pi(g) u_1, u_2) d g d h d [x] \\
=& \int_{X/G} \int_{G \times G} f_1(x g_1) \overline{f_2}(x h)(\pi(h^{-1} g_1) u_1, u_2) d g_1 d h d [x] \\
=& \int_{X/G} \int_{G \times G} f_1(x g_1) \overline{f_2}(x h)(\pi( g_1) u_1, \pi(h) u_2) d g_1 d h d [x]\\
=& \int_{X/G} (\int_G f_1(x g_1) \pi(g_1) u_1 d g_1, \int_G f_2(x h) \pi(h) u_2 d h ) d [x].
\end{split}
\end{equation} 
Here the measure on $X/G$ is normalized so that $d x= d [x] d g$.
Since the second line converges absolutely, by Fubini's Theorem, we can interchange integrals. \\
\\
Observe that the function
$$X \ni x \rightarrow \int_G f_1(x g_1) \pi(g_1) u_1 d g_1$$
is well-defined. We denote it by $\mc I(f_1 \otimes u_1)$. $\mc I(f_1 \otimes u_1)$ is compactly supported and
satisfies that for every $g_2 \in G$,
$$\mc I(f_1 \otimes u_1)(x g_2)=\pi(g_2^{-1}) \mc I(f_1 \otimes u_1)(x) \in \mc H.$$
So $\mc I(f_1 \otimes u_1)$ is a smooth section of the Hilbert vector bundle
$$ X \times_G \mc H \rightarrow X/G$$
with compact support. We obtain a map:
$$\mc I: C_c^{\infty}(X) \otimes \mc H \rightarrow C_c^{\infty}(X \times_G \mc H).$$
By Equation (\ref{gr}), we have 
$$(f_1 \otimes_G u_1)(f_2 \otimes u_2)= (\mc I(f_1 \otimes u_1), \mc I(f_2 \otimes u_2))_{L^2(X \times_G \mc H)}.$$
\\
Let $V$ be a dense subspace of $\mc H$.  We would like to show that 
$${\rm span}\{\mc I(f \otimes u) \mid f \in C_c^{\infty}(X), u \in V \}$$
is dense in $L^2(X \times_G \mc H)$. Suppose otherwise. Then there exists a  $\Psi \in L^2(X \times_G \mc H)$ such that
for any $f \in C_c^{\infty}(X), u \in V$, we have 
$(\mc I(f \otimes u), \Psi)=0$. Notice that 
$$ (\mc I(f \otimes u), \Psi)=\int_{X/G} (\int_G f(xg) \pi(g)u d g, \Psi(x)) d [x]= \int_{X} f(x) (u, \Psi(x)) d x=0.$$
Hence for any fixed $u \in V$, $(u, \Psi(x))=0$ almost everywhere. Let $v \in \mc H$ be an arbitrary vector in $\mc H$. Since $\mc H$ is separable and $V$ is dense in $\mc H$, there exists a sequence $\{ u_i \}_1^{\infty} \subset V$ such that $u_i \rightarrow v$. Now for each $u_i$, $(u_i, \Psi(x))=0$ almost everywhere. Since the set $\{ u_i \}$ is countable, $\{u_i \} \perp \Psi(x)$ almost everywhere. Hence $(v, \Psi(x))=0$ almost everywhere.
Choose an orthonormal  basis $\{v_j \}_{j=1}^{\infty}$ for $\mc H$. For each $v_j$, $(v_j, \Psi(x))=0$ almost everywhere. Hence $(\mc H, \Psi(x))=0$ almost everywhere. It follows that $\Psi(x)=0$ almost everywhere. Therefore,
the set $\{\mc I(f \otimes u) \mid f \in C_c^{\infty}(X), u \in V \}$ spans a dense subspace of $L^2(X \times_G \mc H)$.\\
\\
If $\sum_{i=1}^{l} f_i \otimes_G u_i=0$, by Equation (\ref{gr}), $(\sum_{i=1}^l \mc I(f_i \otimes u_i), \mc I(f \otimes u))=0$ for every $f \in C_c^{\infty}(X)$ and $u \in \mathcal H$. So
$\sum_{i=1}^l \mc I(f_i \otimes u_i)=0$. Conversely, if $\sum_{i=1}^l \mc I(f_i \otimes u_i)=0$,
$\sum _{i=1}^l f_i \otimes_G u_i=0$. We can identify $C_c^{\infty}(X) \otimes_G \mathcal H$ with its image in $C_c^{\infty}(X \times_G \mc H) \subseteq L^2(X \times_G \mathcal H, d x)$. $\Box$ \\

\subsection{Representations obtained from invariant tensor product}
\begin{lem}\label{gkh}
Let $G_1$ be a unimodular group. Let $(\sigma_{1,2}, \mc H_{1,2})$ be a Hilbert representation of $G_1 \times G_2$ and $(\sigma_1, \mc H_1)$  a Hilbert representation of $G_1$. Let $V_{1,2}$ be a $G_2$-stable subspace of $\mc H_{1,2}$ and $V_1$ be a  subspace of $\mc H_1$. Suppose that $V_{1,2}^c$ as a subspace of $ \mc H_{1,2}^*$ is a $\sigma_{1,2}^*(G_2)$-invariant subspace. Suppose that $V_{1,2} \otimes_{G_1} V_1$ is well-defined. Then $V_{1,2} \otimes_{G_1} V_1$ inherits a $G_2$-action from $V_{1,2}$, namely
$$g: v_{1,2} \otimes _{G_1} v_1 \rightarrow \sigma_{1,2}(g) v_{1,2} \otimes_{G_1} v_1, \qquad (v_{1,2} \in V_{1,2}, v_1 \in V_1).$$
Hence $V_{1,2} \otimes_{G_1} V_1$ is  a linear representation of $G_2$. 
\end{lem}
Proof: Consider $\mc A_{G_1}: V_{1,2} \otimes V_1 \rightarrow \Hom (V_{1,2}^c \otimes V_1^c, \mathbb C)$. It suffice to show that
the subspace $ker(\mc A_{G_1})$ is $G_2$-invariant. Suppose that $\sum_{j=1}^t v_{1,2}^j \otimes v_1^j \in ker(\mc A_{G_1})$. Then for any $u_{1,2} \in V_{1,2}^c, u_1 \in V_1^c$, we have
$$\int_{G_1}\sum_{j=1}^t (\sigma_{1,2}(g_1) v_{1,2}^j, u_{1,2})(\sigma_1(g_1) v_1^j, u_1) d g_1=0.$$
Since $V_{1,2}^c$ is a $\sigma_{1,2}^*(G_2)$-invariant subspace,  for any $g_2 \in G_2$ and  $u_{1,2} \in V_{1,2}^c$,  $\sigma_{1,2}^*(g_2^{-1}) u_{1,2} \in V_{1,2}^c \subset \mc H_{1,2}^*$. We have
\begin{equation}
\begin{split}
 & \int_{G_1}\sum_{j=1}^t (\sigma_{1,2}(g_1) \sigma_{1,2}(g_2) v_{1,2}^j, u_{1,2})(\sigma_1(g_1) v_1^j, u_1) d g_1 \\
 = & \int_{G_1}\sum_{j=1}^t(\sigma_{1,2}(g_1)  v_{1,2}^j, \sigma_{1,2}^*(g_2^{-1})u_{1,2})(\sigma_1(g_1) v_1^j, u_1) d g_1=0.
 \end{split}
 \end{equation}
Hence $ \sum_{j=1}^t \sigma_{1,2}(g_2)  v_{1,2}^j \otimes v_1^j \in ker(\mc A_{G_1})$. $\Box$ \\
\begin{re}\label{gkh1}
\begin{enumerate}
\item If $\sigma_{1,2}|_{G_2}$ is unitary and $V_{1,2}$ is $G_2$-invariant, then $V_{1,2}^c$ will automatically be $\sigma_{1,2}^*(G_2)$-invariant. $G_2$-action on $V_{1,2} \otimes_{G_1} V_1$ is compatible with the action of $G_2$ on $\Hom(V_{1,2}^c \otimes V_1^c, \mathbb C)$. 
\item 
For $G_2$ a semisimple Lie group, let $K_2$ be a maximal campact subgroup of $G_2$. The same statement holds for $(\f g_2, K_2)$-module $V_{1,2}$.
\end{enumerate}
\end{re}
Similarly, we have the following two lemma.
\begin{lem}\label{tensorrepn1} Let $ G_2$ be a semisimple Lie group. Let $K_2$ be a maximal compact subgroups of $G_2$. Let $(\sigma_{1,2}, \mc X_{1,2})$ be a continuous representation of $G_1 \times G_2$  and $(\sigma_1, \mc X_1)$ be a continuous representation of $G_1$. Let $V_{1,2}$ be a $(\f g_2, K_2)$-module  in $\mc X_{1,2}$ and $V_1$  a subspace of $\mc X_1$. Let $U_{1,2}$ be  a $(\f g_2, K_2)$-module in $\mc X_{1,2}^*$ and $U_1$ a subspace of $\mc X_1^*$. Then $V_{1,2} \otimes_{G_1, U_{1,2} \otimes U_1} V_1$ inherits a $(\f g_2, K_2)$-module structure from $V_{1,2}$.
\end{lem}
\begin{lem}\label{tensorrepn} Let $ G_2$ be a semisimple Lie group. Let $K_2$ be a maximal compact subgroups of $G_2$. Let $(\sigma_{1,2}, \mc X_{1,2})$ be a smooth representation of $G_1 \times G_2$ equipped with a $(\f g_2, K_2)$ invariant Hermitian form and $(\sigma_1, \mc X_1)$ be a Hermitian representation of $G_1$. Let $V_{1,2}$ be a $(\f g_2, K_2)$-module  in $\mc X_{1,2}$. Let $V_1$ be a subspace of $\mc X_1$. Then $V_{1,2} \otimes_{G_1} V_1$ inherits a $(\f g_2, K_2)$-module structure from $V_{1,2}$.
\end{lem}
Here the requirement that $\mc X_{1,2}$ has a $(\f g_2, K_2)$ invariant Hermitian form will guarantee that  $V_{1,2}^c$ is a $(\f g_2, K_2)$-submodule of $(\sigma_{1,2}^*, \mc X_{1,2}^*)$. Thus, an argument similar to Lemma \ref{gkh} applies.
\begin{ex}[Injectivity]
Let 
$$V_1 \subseteq V_2 \subseteq \ldots \subseteq V_n \subseteq \mathcal X$$
be a sequence of subspaces of a Hilbert representation $(\pi, \mc X)$ of $G$. Let $U$ be a subspace of
the Hilbert representation $(\sigma, \mc Y)$. Suppose that $V_n \otimes_G U$ is well-defined. Then we have
$$V_1 \otimes_{G, V_n^c \otimes U^c} U \subseteq V_2 \otimes_{G, V_n^c \otimes U^c} U \subseteq \ldots \subseteq V_n \otimes_G U.$$
\end{ex}
\subsection{Invariant Hermitian Forms for ITP} 
Suppose that the representations $(\sigma_{1,2}, \mc X_{1,2})$ and $(\sigma_1, \mc X_1)$ are equipped with $G_1$-invariant Hermitian forms $(\, , \,)$. Notice that
\begin{equation}\label{hermitian}
\begin{split}
(\sum_{j=1}^s v_{1,2}^j \otimes_{G_1} v_1^j)(\sum_{i=1}^t u_{1,2}^i \otimes u_1^i)= & \sum_{i,j} \int_{G_1}(\sigma_{1,2}(g_1) v_{1,2}^j, u_{1,2}^i)(\sigma_1(g_1) v_1^j, u_1^i) d g_1 \\
=& \sum_{i,j} \int_{G_1}( v_{1,2}^j, \sigma_{1,2}(g_1^{-1})u_{1,2}^i)( v_1^j, \sigma_1(g_1^{-1}) u_1^i) d g_1 \\
= & \sum_{i,j} \int_{G_1}( v_{1,2}^j, \sigma_{1,2}(g_1)u_{1,2}^i)( v_1^j, \sigma_1(g_1) u_1^i) d g_1 \\
= & \sum_{i,j} \int_{G_1} \overline{( \sigma_{1,2}(g_1)u_{1,2}^i,  v_{1,2}^j)} \overline{(  \sigma_1(g_1) u_1^i, v_1^j)} d g_1 \\
= & \overline{(\sum_{i=1}^t u_{1,2}^i \otimes_{G_1} u_1^i )(\sum_{j=1}^s v_{1,2}^j \otimes v_1^j)}.
\end{split}
\end{equation}
Here $v_{1,2}^j \in  V_{1,2}$, $v_1^j \in V_1$, $u_{1,2}^i \in V_{1,2}^c$ and $u_1^i \in V_1^c$.
By Lemma \ref{pairing},  $(\sum_{j=1}^s v_{1,2}^j \otimes_{G_1} v_1^j)(\sum_{i=1}^t u_{1,2}^i \otimes u_1^i)$ yields a pairing between $\sum_{j=1}^s v_{1,2}^j \otimes_{G_1} v_1^j$ and $\sum_{i=1}^t u_{1,2}^i \otimes_{G_1} u_1^i$. By our computation, this pairing produces a Hermitian form on $V_{1,2} \otimes_{G_1} V_1$. We give the following

\begin{defn}\label{hermitianform} Suppose that the representations $(\sigma_{1,2}, \mc X_{1,2})$ and $(\sigma_1, \mc X_1)$ are equipped with $G_1$-invariant Hermitian forms $(\, , \,)$. Define a Hermitian form on $V_{1,2} \otimes_{G_1} V_1$ as follows
$$(v_{1,2} \otimes_{G_1} v_1, u_{1,2} \otimes_{G_1} u_1)=(v_{1,2} \otimes_{G_1} v)(u_{1,2} \otimes u_1)=
\int_{G_1}(\sigma_{1,2}(g_1) v_{1,2}, u_{1,2})(\sigma_1(g_1) v_1, u_1) d g_1.$$
\end{defn}

\begin{Proposition}\label{invher} Let $G_1, G_2$ be two semisimple Lie groups. Let $K_i$ be a maximal compact subgroup of $G_i$ respectively. Let $(\sigma_{1,2}, \mc X_{1,2})$ be a continuous representation of $G_1 \times G_2$ equipped with a $G_1 \times G_2$-invariant Hermitian form $(\, , \,)$ and $(\pi_1, \mc X_1)$ be a continuous representation of $G_1$ equipped with a $G_1$-invariant Hermitian form $(\, , \,)$. Let $V_{1,2}$ be a $(\f g_2, K_2)$-module  in $\mc X_{1,2}$. Let $V_1$ be a subspace of $\mc X_1$. Then the Hermitian form on $V_{1,2} \otimes_{G_1} V_1$ is non-degenerate, invariant under the action of both $K_2$ and $\f g_2$.
\end{Proposition}
Proof: Let $u_{1,2}, v_{1,2} \in V_{1,2}$ and $u_1, v_1 \in V_1$. Let $k_2 \in K_2$. Then
\begin{equation}
\begin{split}
 & (\sigma_{1,2}(k_2) u_{1,2} \otimes_{G_1} u_1, \sigma_{1,2}(k_2) v_{1,2} \otimes_{G_1} v_1) \\
 =& \int_{G_1} (\sigma_{1,2}(g_1 k_2) u_{1,2}, \sigma_{1,2}(k_2) v_{1,2})(\sigma_1(g_1) u_1, v_1) d g_1\\
 = &
\int_{G_1} (\sigma_{1,2}(g_1 ) u_{1,2}, v_{1,2})(\sigma_1(g_1) u_1, v_1) d g_1\\
=&(u_{1,2} \otimes_{G_1} u_1,  v_{1,2} \otimes_{G_1} v_1)
\end{split}
\end{equation}
The invariance of $\f g_2$-action is similar. $\Box$.
\begin{cor} Under the hypothesis of Prop. \ref{invher}, the averaging operator
$$\mc A_{G_1}: V_{1,2} \otimes V_1 \rightarrow \Hom (V_{1,2}^c \otimes V_1^c, \mathbb C)$$
is a $(\f g_2, K_2)$-module homomorphism. 
\end{cor}

\section{Representation Theory of $SO(p,q)$}
From now on all representations are assumed to be continuous.
Let $G$ be a semisimple Lie group. Fix a maximal compact subgroup $K$. {\bf A Hilbert representation $(\pi, \mc H)$ of $G$  is said to be admissible if it restricts to a unitary representation of $K$ and the isotypic subspaces $\mc H_{\sigma}$ where $\sigma \in \h{K}$ are finite dimensional}. Let $\mc H_K$ be the $K$-finite vectors of $\mc H$. Then $\mc H_K$ is a $(\f g, K)$-module satisfying certain compatibility conditions. We call $\mc H_K$ the Harish-Chandra module of $\mc H$. We may also call a continuous representation admissible, if (1) the $K$-finite vectors are dense and smooth in the representation, (2) the $K$-finite vectors yield a $(\f g, K)$-module and (3) each $K$-type has finite multiplicity. In all cases, $V_K$ or $V(\pi)$ will be used to denote the $K$-finite vectors of a representation $(\pi, V)$ of $G$. Two admissible representations are said to be infinitesimally equivalent if their underlying Harish-Chandra modules are equivalent. In this paper, we will mainly be interested in infinitesimal equivalences.\\
\\
Let $U(\f g)$ be the universal enveloping algebra. Let $\mc Z$ be the center of $U(\f g)$. Let $I: \mc Z \rightarrow \mathbb C$ be a character of $\mc Z$. We say that a $(\f g, K)$-module $V$ has infinitesimal character $I$ if $\mc Z$ acts on $V$ by the character $I$. We use the Harish-Chandra isomorphism to identify the spectrum of $\mc Z$ with the complex dual of $\f h$, a maximal Cartan subalgebra of $\f g_{\mathbb C}$. This identification is unique up to the action of the Weyl group $W(\f g_{\mathbb C}, \f h)$. In other words, $I$ and $w I$ are the same infinitesimal character if
$w \in W(\f g_{\mb C}, \f h)$. If $(\pi, V)$ is an admissible representation of $G$ and $V_K$ has infinitesimal character $\lambda$, then we say that $\pi$ has infinitesimal character $\lambda$.\\
\\
{\bf Let $G=SO_0(p,q)$ be the identity component of $SO(p,q)$}. Suppose that $pq \neq 0$. So $G$ is  noncompact. We often assume $q \geq p$ so the noncompact component of the $KAK$-decomposition is in $\mathbb R^p$. This is not necessary, but convenient. If $q < p$, then the noncompact component of the $KAK$ decomposition is in $\mathbb R^q$. \\
\\
{\bf Let $\chi$ be the unique  unitary character of $SO(p,q)$ that maps the nonidentity component to $-1$. We extend $\chi$ to a character of products of real general linear groups and $SO(p,q)$ by defining
$$\chi|_{GL}=\frac{\det}{| \det |}.$$}
Define $\chi$ on $SO(k,0)$ to be the trivial character. Take $\chi$ to be the universal character of $SO$, $GL$ and their products.

\subsection{Representations of SO(p,q)}
Let $(\pi, V)$ be an irreducible admissible representation of $G$. Let $h$ be an element in the nonidentity component of $SO(p,q)$ such that $h^2=I$. Then $SO(p,q)$ can be identified with a semidirect product of $\{ I, h \}$ and $G$. Define
an irreducible admissible representation $(\pi^h, V)$ as
$$\pi^h(g)=\pi(h^{-1} g h) \qquad (g \in G).$$
Then $\pi^h \cong \pi$ or $\pi^h \ncong \pi$.
\subsubsection{The case $\pi^h \cong \pi$}
If $\pi^h \cong \pi$, then there is an intertwining operator $i: V \rightarrow V$ such that
$$\pi(g) i= i \pi^h(g)=i \pi(h^{-1} g h), \qquad (g \in G).$$
It follows that $$\pi(g) i i= i \pi(h^{-1} g h) i= i i \pi(h^{-2}g h^2)=i i \pi(g).$$
Since $\pi$ is irreducible, by Schur's Lemma, $i i = \lambda \in \mathbb C-\{0\}$. Normalize $i$ so that $i^2=I$.
Define a representation $(\overline{\pi}, V)$ of $SO(p,q)$ such that
$$\overline{\pi}(h)=i, \qquad \overline{\pi}(g)=\pi(g) \, \qquad (\forall \, g \in G).$$
It is easy to check that $\overline{\pi}$ is an admissible irreducible representation of $SO(p,q)$. $\overline{\pi}$ has the following properties:
\begin{enumerate}
\item $\overline{\pi}|_G \cong \pi$;
\item $\overline{\pi} \otimes \chi \ncong \overline{\pi}$.
\item $\overline{\pi}$ and $\overline{\pi} \otimes \chi$ come from
two different choices of normalization.
\end{enumerate}
\subsubsection{The case $\pi^h \ncong \pi$}
If $\pi^h \ncong \pi$, then we can define an admissible  representation
$(\overline{\pi}, V \oplus V)$ as follows:
$$\overline{\pi}(h)(v_1, v_2)=(v_2, v_1), \qquad \overline{\pi}(g)(v_1, v_2)=(\pi(g) v_1, \pi(h^{-1} g h) v_2) \qquad (v_1, v_2 \in V, g \in G).$$
Then
$$\overline{\pi}(h^{-1})\overline{\pi}(g)\overline{\pi}( h)(v_1, v_2)=\overline{\pi}(h)\overline{\pi}(g)(v_2, v_1)=\overline{\pi}(h)(\pi(g)v_2, \pi(h^{-1}g h) v_1)$$
$$ =(\pi(h^{-1} g h) v_1, \pi(g) v_2)=\overline{\pi}(h^{-1} g h)(v_1, v_2).$$
Hence $\overline{\pi}(h^{-1})\overline{\pi}(g)\overline{\pi}( h)=\overline{\pi}(h^{-1} g h)$. Obviously, $\overline{\pi}(h)^2=I$.
So $\overline{\pi}$ is an admissible representation of $SO(p,q)$.
It is irreducible by Mackey analysis. $\overline{\pi}$ has the following properties:
\begin{enumerate}
\item $\overline{\pi}|_G= \pi \oplus \pi^h $ and $\pi \ncong \pi^h$;
\item $\overline{\pi} \otimes \chi \cong \overline{\pi}$.
\end{enumerate}
Similar statements  hold for irreducible unitary  representations of the compact group $O(m)$ and the general linear group $GL(n)$. Unless otherwise specified, if $\pi$ is an irreducible admissible representation of $G$, $\overline{\pi}$ will be used to denote the one, or one of the two representations defined in this section. 

\subsection{ Representations of $SO(p)$ and $O(p)$}\label{compactrepn}
Given a compact Lie group $K$, let $\h{K}$ be the the equivalence classes of irreducible unitary representations of $K$. When $K$ is connected, $\hat{K}$ is parametrized by integral dominant weights upon a choice of a positive root system. If $\xi$ is the highest weight of $\pi$, we write $(\pi_{\xi}, V_{\xi})$ or sometimes just $(\xi)$ for  $\pi$. \\
\\
Let $K=SO(p) \times SO(q)$. The irreducible unitary representations of $K$ are products of irreducible unitary representations of $SO(p)$ and $SO(q)$.
\subsubsection{Odd Orthogonal Groups}
For $p$ odd, $\h{SO}(p)$ is parameterized by an integral vector $\xi$ of dimension $[\frac{p}{2}]$:
$$\xi_1 \geq \xi_2 \geq \ldots \geq \xi_{[\frac{p}{2}]-1} \geq \xi_{[\frac{p}{2}]} \geq 0.$$
Denote this set of integral vectors  by $\Pi_p(\mathbb Z)^+$.
In this case, $-I_p \in O(p)-SO(p)$ and $O(p) \cong SO(p) \times \{ \pm I_p \}$. The unitary dual $ \h{O}(p)$ can be parametrized by $(\xi, \pm)$ with $\xi \in \Pi_p(\mathbb Z)^+$. Here $\pi_{\xi, \pm}|_{SO(p)}$ is $\pi_{\xi}$ and $\pi_{\xi, \pm}(-I_p)$ is $\pm$ identity.  \\
\\
When $p=1$, there is only a trivial representation of $SO(1)$. We denote it by $(0)$. There are two irreducible unitary representations of $O(1)$. We denote them by $(0, +)$ and $(0, -)$. \\
\subsubsection{Even Orthogonal Group}
For $p$ even, $\h{SO}(p)$ is parametrized by an integral vector
$\xi$ of dimension $\frac{p}{2}$:
$$\xi_1 \geq \xi_2 \geq \ldots \geq \xi_{\frac{p}{2}-1} \geq |\xi_{\frac{p}{2}}|.$$
Denote this set of integral vectors by $\Pi_p(\mathbb Z)^+$.
Fix an element $h \in O(p)-SO(p)$ such that $h^2=I_p$.
\begin{enumerate}
\item
If $\xi_{\frac{p}{2}} \neq 0$, then $\pi^{h}_{\xi} \cong \pi_{\xi^{-}}$ with
$\xi^{-}=(\xi_1, \xi_2, \ldots \xi_{\frac{p}{2}-1}, -\xi_{\frac{p}{2}}).$
We obtain an irreducible unitary representation of $O(p)$ such that
$\overline{\pi}|_{SO(p)}=\pi \oplus \pi^h.$
Following \cite{kv}, we parametrize such $\overline{\pi}$ by {\bf $(|\xi|, +)$}.
\item
If $\xi_{\frac{p}{2}} =0$, then $\pi^h_{\xi} \cong \pi_{\xi}$. So 
$\pi$ extends to two representations of $O(p)$, differing by $\det$. we shall follow the convention set in (6.10 \cite{kv}) and parametrize these two representations by ${(\xi, +)}$ and ${(\xi, -)}$. For example, the standard representation of $O(p)$ will be parametrized by
$$(\xi, +), \qquad \xi=(1, 0, \ldots, 0).$$
\end{enumerate}
The distinction between $(\xi, +)$ and $(\xi, -)$ will become crucial when we apply Howe's theory of dual pair correspondence.  
\subsubsection{Representations of $S(O(p)O(q))$}
In our notation, we will have two kinds of irreducible representations of $O(p)$. One kind is denoted by $(\xi, +)$ with
$\xi_1 \geq \xi_2 \ldots \geq \xi_{[\frac{p}{2}]} \geq 0$. This kind will include 2.2.2 (1). The other kind is denoted by $(\xi, -)$ with $\xi_1 \geq \xi_2 \ldots \geq \xi_{[\frac{p}{2}]} \geq 0$.\\
\\
Let $S(O(p)O(q))=\{(k_1, k_2) \in O(p) O(q) \mid \det k_1 \det k_2=1 \}$ with $pq \neq 0$. $S(O(p)O(q))$ is a maximal compact subgroup of $SO(p,q)$. Let $\sigma$ be an irreducible representation of $S(O(p)O(q))$. Let $(\xi) \otimes (\eta)$ be an irreducible subrepresentation of $\sigma|_{SO(p)SO(q)}$. \\
\\
If $(\xi)$ and $(\eta)$ both extend to  irreducible representations of $O(p)$ and $O(q)$, then $$(\xi, +) \otimes (\eta, +)|_{S(O(p)O(q))} \cong (\xi, -) \otimes (\eta, -)|_{S(O(p)O(q))}.$$ We parametrize this irreducible representation of $S(O(p)O(q))$ by
$(\xi, \eta, +)$. We parametrize $(\xi, +) \otimes (\eta, -)|_{S(O(p)O(q))} \cong (\xi, -) \otimes (\eta, +)|_{S(O(p)O(q))}$ by
$(\xi, \eta, -)$. \\
\\
Otherwise, $\sigma|_{SO(p)SO(q)}$ will have two irreducible constituents, with one being $(\xi) \otimes (\eta)$. we parametrize $\sigma$ by $(|\xi|, |\eta|, +)$.

\subsection{Basic Structure Theory of $G$}
Suppose now that $p \leq q$. Fix a maximal compact subgroup $K=SO(p)SO(q)$ in $G$. Fix an Iwasawa decomposition $KA_{\rm m}N_{\rm m}$.
Let $\Delta(\f g, \f a_{\rm m})$ be the restricted roots and $\Delta^+(\f g, \f a_{\rm m})$ be the positive restricted roots defined by $N_{\rm m}$. Let $\{e_1, e_2, \ldots, e_p \}$ be the standard basis for $\f a_{\rm m}^*$ such that
$$\Delta(\f g, \f a_{\rm m})= \{ \pm e_i \pm e_j \mid i > j,  i, j \in [1, p] \} \qquad (p=q);$$
$$\Delta(\f g, \f a_{\rm m})=\{ \pm e_i \pm e_j \mid i > j, i, j \in [1, p] \} \cup \{\pm e_i \mid i \in [1, p] \}^{q-p} \qquad (q > p).$$
Here $q-p$ denotes the multiplicities. Let $\rho(p,q)$ be the half sum of positive restricted roots. Then
$$\rho(p,q)=( \frac{q+p}{2}-1, \frac{q+p}{2}-2, \ldots, \frac{q-p}{2} ).$$
{\bf The basis $\{ e_i \}_{i=1}^p$ defines a coordinate system for $\f a_{\rm m}$. By this coordinate system, we identify $\f a_{\rm m}$ with $\mathbb R^p$. Equip $\f a_{\rm m}$ and $\f a_{\rm m}^{*}$ with the standard inner product of $\mathbb R^p$.} The open Weyl chamber is
$$\f a_{\rm m}^{+}=\{ H_1 > H_2 > \ldots > H_p > 0 \mid H \in \mathbb R^p \} \qquad (q> p);$$
$$\f a_{\rm m}^{+}=\{ H_1 > H_2 > \ldots > | H_p|  \mid H \in \mathbb R^p \} \qquad (q=p).$$
For $q >p$, the Weyl group $W(\f g,\f a_{\rm m})$ is generated by permutations and sign changes. For $q=p$, the Weyl group $W(\f g, \f a_{\rm m})$ is generated by permutations and even number of sign changes. Let $\overline{A_{\rm m}}^+= \exp \overline{\f a_{\rm m}}^+$.  \\
\\
Consider now the Cartan decomposition $KA_{\rm m}K$. Since $W(\f g, \f a_{\rm m})=N_K(\f a_{\rm m})/Z_K(\f a_{\rm m})$, $G$ can be decomposed
as $K \overline{A_{\rm m}}^+ K$. Here $N_K(\f a_{\rm m})$ is the normailizer of $\f a_{\rm m}$ in $K$ and $Z_K(\f a_{\rm m})$ is the centralizer of $\f a_{\rm m}$ in $K$. For any $g \in G$, write $g=k_1 \exp H^{+}(g) k_2$ with $H^+ \in \overline{\f a_{\rm m}}^+$. $H^+(g)$ is uniquely determined by $ g$. \\
\\
Let $M_m=Z_K(\f a_{\rm m})$. Then $M_m$ can be identified with 
$$\{ (\epsilon_1, \epsilon_2, \ldots \epsilon_p, h) \mid \epsilon_i=\pm 1, h \in SO(q-p), \prod_{i=1}^{p} \epsilon_i =1 \}.$$
Let $P_{\rm m}=M_{\rm m}A_{\rm m}N_{\rm m}$. $P_{\rm m}$ is a minimal parabolic subgroup.

\subsection{Induced Representations }
Let $P \supseteq P_{\rm m}$ be a parabolic subgroup of $G$. Let $P=LN$ be the Levi decomposition. Then $L$ has the following form
$$\{ \chi (g_0) \prod_{i=1}^l \chi (h_i)=1 \mid ( h_1, h_2, \ldots h_l, g_0) \in GL(r_1)GL(r_2) \ldots GL(r_l) SO(p-\sum_{i=1}^l r_i, q-\sum_{i=1}^l r_i) \}.$$
where $r_i \geq 1 (i \in [1, l]),  p \geq \sum_{i=1}^l r_i$. Let $P=MAN$ be the Langlands decomposition. {\bf Let $SL^{\pm}(r_i)=\{\det g = \pm 1 \mid g \in GL(r_i) \}.$}
 Then
 $A \cong ({\mathbb R}^+)^l$ and 
 $$M= \{ \chi ( g_0) \prod_{i=1}^l \det h_i=1 \mid ( h_1, h_2, \ldots h_l, g_0) \in \prod_{i=1}^l SL^{\pm}(r_i) \times SO(p-\sum_{i=1}^l r_i, q-\sum_{i=1}^l r_i)\}.$$
 The connected component of $M$, 
 $M_0 \cong SL(r_1) SL(r_2) \ldots SL(r_l) SO_0(p-\sum_{i=1}^l r_i, q-\sum_{i=1}^l r_i)$.\\
 \\
 For the special orthogonal group $SO(p,q)$, the Levi factor of a parabolic subgroup must be of the form 
 $$GL(r_1)GL(r_2) \ldots GL(r_l) SO(p-\sum_{i=1}^l r_i, q-\sum_{i=1}^l r_i);$$
 and $M \cong \prod_{i=1}^l SL^{\pm}(r_i) \times SO(p-\sum_{i=1}^l r_i, q-\sum_{i=1}^l r_i)$. \\
 \\
 Let $\Delta^{+}(\f g, \f a)$ be the positive roots from $N$. Let $\rho$, or sometimes $\rho_{P}$ be the half sum of  the positive roots (with multiplicities) in $\Delta^+(\f g, \f a)$. Notice that $\Delta(\f g, \f a)$ in general is not a root system. Nevertheless, $\Delta^{+}(\f g, \f a)$ still defines an open positive Weyl chamber in $\f a^{*}$, namely
 $$[{\f a}^{*}]^+=\{ v \in \f a^* \mid  (v, \alpha) > 0, \,\, \forall \, \alpha \in \Delta^+(\f g, \f a) \,\} .$$
 Let $(\sigma, \mc H_{\sigma})$ be an admissible Hilbert representation of $M$. Let $v \in (\f a^*)_{\mathbb C}$. Let
$\Ind_{P}^G \sigma \otimes \mathbb C_v$ be the normalized induced representation. The Hilbert space of $\Ind_{P}^G \sigma \otimes \mathbb C_v$ is given by
$$\{ f: G \rightarrow \mc H_{\sigma} \mid f(g m a n)=\exp -(v+\rho)(\log a) \sigma(m)^{-1} f(g), f|_K \in L^2(K, \mc H_{\sigma}) \}.$$
If $(\sigma, V_{\sigma})$ is an admissible representation of $M$, we can also take continuous functions 
$$\{ f: G \rightarrow V_{\sigma} \, continuous \mid f(g m a n)=\exp -(v+\rho)(\log a) \sigma(m)^{-1} f(g)  \}.$$
Again, we obtain an admissible representation of $G$. Whether taking $L^2$-sections or continuous sections, the underlying Harish-Chandra module is the same, namely, $K$-finite vectors in
$\Ind_{M \cap K}^K [V_{\sigma}]_{M \cap K} .$
So the two induced representations are infinitesimally equivalent. \\
\\
Parabolically induced representations have the following nice properties:
\begin{enumerate}
\item If $\Re v=0$ and $\sigma$ is unitary, the $\Ind_P^G \sigma \otimes \mathbb C_{v}$ is unitary. This type of induction is called unitary parabolic induction.
\item If $\sigma$ has infinitesimal character $\lambda$, then $\Ind_{P}^G \sigma \otimes \mathbb C_v$ has infinitesimal character $\lambda \oplus v$.
\item If $\sigma|_{M \cap K}$ has a composition series of finite length, then $(\Ind_{P}^G \sigma \otimes \mathbb C_v)_K$ will also have a composition series of finite length.
\item $\Ind_P^G \sigma \otimes \mathbb C_v|_K=\Ind_{M \cap K}^K (V_{\sigma})_{M \cap K}$. The action of $K$ is independent of $v$.
\item Every irreducible admissible representation occurs infinitesimally as a subrepresentation of a principal series representation, namely those induced representation with $P=P_{\rm m}$ and $\sigma$ irreducible.
\end{enumerate}
\subsection{Langlands Classification }
We shall now review some basic facts about Langland classification. Langlands reduces the classification of the infinitesimal equivalence classes of irreducible admissible representations to the classification of irreducible tempered representations (\cite{la} \cite{kn}). Then Knapp and Zuckerman complete the classification of irreducible tempered representations, which are necessarily unitarizable. 
\begin{Theorem}[Langlands Classification]
Fix a minimal parabolic subgroup $P_{\rm m}$. Then the infinitesimal equivalence classes of irreducible admissible representations of $G$ is in one-to-one correspondence to the following triples
\begin{enumerate}
\item $P=MAN$, a parabolic subgroup containing $P_{\rm m}$;
\item the equivalence class of $\sigma$, where $\sigma$ is an irreducible tempered representation of $M$;
\item $v \in (\f a^*)_{\mathbb C}$ such that $\Re(v) \in [\f a^*]^+$.
\end{enumerate}
The irreducible admissible representation $J(P,\sigma, v)$ is given by the unique irreducible admissible quotient of $\Ind_{MAN}^{G} \sigma \otimes \mathbb C_{v}$.
\end{Theorem}
In particular, Langlands classification says that there is a unique irreducible quotient for $\Ind_{MAN}^{G} \sigma \otimes \mathbb C_{v}$. So $J(P, \sigma, v)$ is often called the Langlands quotient. Notice that if $\Re (v) \notin [\f a^*]^+$ or $\sigma$ not tempered, there could be many irreducible quotients. At the end of this section, we will \lq\lq isolate \!\rq\rq some Langlands quotient without assuming that $\Re(v)$ is in the open Weyl chamber, using Vogan's lowest $K$-types. 
\begin{re}
There is a revised Langlands classification based on cuspidal parabolic subgroups. I shall only state some facts that will be used later. Let $P=MAN$ be a cuspidal parabolic subgroup containing $P_{\rm m}$. Let $\sigma$ be a discrete series representation of $M$. Suppose that $\Re(v)$ is in the {\bf closed} positive Weyl chamber. Then there is a Langlands quotient $J(P, \sigma, v)$, defined to be the image of certain intertwining operator.  In this case, $J(P, \sigma, v)$ may not be irreducible. Nevertheless, it decomposes into a direct sum of irreducibles.  Every irreducible Harish-Chandra module can be constructed as a direct summand of $J(P, \sigma, v)$. 
\end{re}
\subsection{Matrix Coefficients }
Now we shall list a few important properties, that more or less characterize Langlands quotients. We will start with irreducible tempered representations. These representations are all unitarizable. {\bf  Let $\Xi_{\lambda}(g)$ be Harish-Chandra's spherical function. We write
$\Xi(g)$ for $\Xi_0(g)$}. Notice that $\Xi_{\lambda}(g)=\Xi_{w \lambda}(g)$ for any $w \in W(\f g, \f a_{\rm m})$.
\begin{Theorem} Let $\pi$ be an irreducible admissible representation of $G$.
The following are equivalent.
\begin{enumerate}
\item $(\pi, \mc H)$ is an irreducible tempered representation.
\item All $K$-finite matrix coefficients of $\pi$ are bounded by a multiple of $\Xi(g)$.
\item All $K$-finite matrix coefficients of $\pi$ are in $L^{2+\epsilon}(G)$ for every $\epsilon >0$.
\item $(\pi, \mc H)$ is infinitesimally equivalent to a subrepresentation of $\Ind_P^G \sigma \otimes \mathbb C_v$ with $\sigma$ a discrete series and $\Re(v)=0$.
\item All leading exponent $v$ of $\pi$ satisfies $\Re(v+\rho(p,q))(\f a_{\rm m}^+) \leq 0$.
\end{enumerate}
\end{Theorem}
If $(\pi, \mc H)$ is assumed to be unitary, then the conditions above are equivalent to the condition that $(\pi, \mc H)$ is weakly contained in $L^2(G)$ (\cite{chh} \cite{wallach}).
Now we make some remarks about Langlands quotient.
\begin{re}\label{langlandsmatrixco}
\begin{enumerate}
\item The $K$-finite matrix coefficients of $J(P, \sigma, v)$ are bounded by multiples of $\Xi_{\Re(v) \oplus 0}(g)$. Here we embed $\f a^*  \rightarrow \f a_{\rm m}^* $.
\item The $K$-finite matrix coefficients of $\Ind_P^G \sigma \otimes \mathbb C_v$ are bounded by multiples of $\Xi_{\Re(v) \oplus 0}(g)$ (see Prop 7.14 \cite{kn}).
\item If $\Re(v) \in -[\f a^*]^+$, then there is a unique irreducible admissible subrepresentation in $\Ind_P^G \sigma \otimes \mathbb C_v$. This representation is the dual representation of $J(P, \sigma, -{v})$. 
\end{enumerate}
\end{re}
Langlands' classification can be established by studying the leading exponents of the asymptotic expansion of matrix coefficients. We quote some result from \cite{kn}.
\begin{Theorem}\label{leadingex} Let $\pi$ be an irreducible admissible representation of $G$. Let $v-\rho(p,q)$ be a leading exponent of $\pi$. Let $\f t_m$ be a maximal Cartan subalgebra of $\f o(q-p)$. Then $v$ is related to the infinitesimal character $\mc I(\pi) \in (\f a_{\rm m} \oplus \f t_m)^*_{\mathbb C}$ by the following equation
$$v = w \mc I(\pi)|_{\f a_{\rm m}} $$ 
where $w$ is an element in the Weyl group of the complex Lie algebra.  
If there is a $\lambda$ such that each leading exponent $v-\rho(p,q)$ satisfies the condition that
$$ | \rho(p,q)-v| \succeq -\lambda $$
then every $K$-finite matrix coefficient $f(k_1 \exp H^+ k_2)$ is bounded by a multiple of 
$$(1+(H^+, H^+))^Q \exp(\lambda(|H^+|))$$
for some $Q >0$.
\end{Theorem}

\subsection{Vogan's Subquotient}
In practice, Langlands classification tells little about the algebraic structure of the Harish-Chandra module. It is not easy to apply Langlands classification, for example, to determine the composition factors of an admissible quasisimple representation. In \cite{vogan79}, Vogan took a more algebraic approach and studied the $K$-types of an irreducible admissible representation.  Vogan proved that in each irreducible admissible representation each minimal $K$-type appears with multiplicity one. 

\begin{defn} Let $\pi$ be an admissible representation of $G$. Fix a maximal compact subgroup $K$. Let $V(\pi)$ be the Harish-Chandra module of $\pi$. Suppose that $\mu$ is a lowest $K$-type with multiplicity one in $V(\pi)$. Let $V_0(\pi, \mu)$ be the $(\f g, K)$-module generated by $\mu$ in $V(\pi)$. Let $V_1(\pi, \mu)$ be the maximal $(\f g, K)$-submodule of $V_0(\pi, \mu)$ not containing $\mu$. We call $V_0(\pi, \mu)/V_1(\pi, \mu)$ a Vogan subquotient of $V(\pi)$. If Vogan subquotient is unique in $V(\pi)$, we say that $V_0(\pi, \mu)/V_1(\pi, \mu)$ is the Vogan subquotient of $V(\pi)$.
\end{defn}
We shall make several remarks here.
\begin{enumerate}
\item
First, by definition,  a Vogan subquotient is an irreducible $(\f g, K)$-module. 
\item The same definition can be generalized to allow lowest $K$-types with multiplicities. Then $V_0(\pi, \mu)/V_1(\pi, \mu)$ may not be irreducible.
\item The lowest $K$-type may not be unique, even for irreducible Harish-Chandra modules.
\end{enumerate}

\begin{Theorem}[\cite{vogan79}] Let $P=MAN$ be a cuspidal parabolic subgroup of $G$ containing the minimal parabolic subgroup $P_{\rm m}$. Let $\sigma$ be a discrete series representation of $G$. Suppose that $\Re(v)$ is in the closed positive Weyl chamber. Then $J(P, \sigma, v)$ is equivalent to the direct sum of Vogan subquotients of $\Ind_P^G  \sigma \otimes \mathbb C_{v}$.
\end{Theorem}

\begin{defn} Let $G$ be a real reductive group.
Given a cuspidal parabolic subgroup $P$, an irreducible  tempered representation $\sigma$, any $v \in \f a^*_{\mathbb C}$,
if  $\Ind_P^G  \sigma \otimes \mathbb C_{v}$ has a unique Vogan subquotient $\pi$, we call $(P, \sigma, v)$ the Langlands-Vogan parameter of $\pi$. Here $v$ is not unique. 
\end{defn}
\section{Degenerate Principal Series $I_n(s)$}
Let $s \in \mathbb R$. Consider $G=SO_0(n,n)$. Let $I_n(s)=\Ind_{GL(n)_0 N}^{SO_0(n,n)} |\det|^s$. Notice that
$$[\Ind_{GL(n) N}^{SO(n,n)} |\det|^s \otimes \mathbb C_{\pm}] |_{SO_0(n,n)} \cong \Ind_{GL(n)_0 N}^{SO_0(n,n)} |\det|^s,$$
Here $\mathbb C_{+}$ is the trivial representation of $SL^{\pm}(n)$ and $\mathbb C_{-}$ is the sign character  of $SL^{\pm}(n)$.   In addition
$$\Ind_{GL(n) N}^{SO(n,n)} |\det|^s \otimes \mathbb C_{-} \cong \chi \otimes [\Ind_{GL(n) N}^{SO(n,n)} |\det|^s \otimes \mathbb C_{+}].$$
Essentially, for Siegel parabolic subgroups, there is only one degenerate principal series of $SO(n,n)$ to study. If one considers the group $Spin_{0}(n,n)$,
there is another degenerate principal series. The composition factors, $K$-types and unitarity of the composition factors are determined by K. Johnson (\cite{jo}) and reworked by several authors in a more general context (\cite{sa} \cite{zhang} \cite{ll}). {\bf In this section, our quotients and subrepresentations shall be interpreted in the category of Harish-Chandra modules.}
\subsection{Reducibility, $K$-types, composition series and Unitarizability}
Recall that the infinitesimal character of $I_n(s)$,
$$\mc I(I_n(s))=(\frac{n-1}{2}+s, \frac{n-3}{2}+s, \ldots, -\frac{n-3}{2}+s, -\frac{n-1}{2}+s).$$
Here the infinitesimal character is unique only up to the action of Weyl group. Due to the nondegenerate pairing between $I_n(s)$ and $I_n(-s)$, we have
$I_n(s) \cong \mc I_n(-s)^*$. There is  an intertwining operator between $I_n(s)$ and $I_n(-s)$. The composition factors of $I_n(s)$ and $I_n(-s)$ are the same. As shown in \cite{jo}, $I_n(s)$ is irreducible if $s \notin \mathbb Z+\frac{n-1}{2}$.
\begin{Theorem}[Johnson, see also \cite{sa} \cite{zhang} \cite{ll}]\label{basicins}
Suppose that $s \in \frac{n-1}{2}+ \mathbb Z$.
\begin{enumerate}
\item The $K$-types of $I_n(s)$, independent of $s$, are multiplicity free and are parameterized by $\{(\lambda, \lambda) \mid \lambda \in \Pi_{n}(\mathbb Z)^+ \}$. More precisely, $K$-types in $I_n(s)$ are of the following pairs of integral highest weights
$$(\lambda_1 \geq \lambda_2 \geq \ldots \geq \lambda_{n_0}) \geq 0, (\lambda_1 \geq \lambda_2 \geq \ldots \geq \lambda_{n_0}) \geq 0, \qquad (K=SO(2n_0+1) SO(2n_0+1));$$
$$(\lambda_1 \geq \lambda_2 \geq \ldots \geq |\lambda_{n_0}|), (\lambda_1 \geq \lambda_2 \geq \ldots \geq |\lambda_{n_0}|), \qquad (K=SO(2n_0)SO( 2n_0)).$$
This is essentially the Peter-Weyl Theorem.
\item For  $n $ odd, $I_n(s)$ has at most $\frac{n-1}{2}+1$ composition factors, namely, $V_i(s) (i \in [0, \frac{n-1}{2}])$. If  $i \geq \frac{n-1}{2}-|s|$, the composition factor $V_i(s)$  exists. It has the following $K$-types
$$\{ (\lambda, \lambda) \mid \lambda_i \geq |s|-\frac{n-1}{2}+i \geq \lambda_{i+1} \}.$$
If $i \neq \frac{n-1}{2}$, $V_i(s)$ is called a small constituent. $V_{\frac{n-1}{2}}(s)$ is called a large constituent.
\item For  $n $ even, $I_n(s)$ has at most $\frac{n}{2}+2$ composition factors, namely, 
$$ V_0(s), v_1(s), \ldots, V_{\frac{n}{2}-2}(s), V_{\frac{n}{2}-1}(s), V_{\frac{n}{2}}(s)^{\pm}.$$
 \begin{enumerate}
 \item Let $i \in [0, \frac{n}{2}-1]$. If  $i \geq \frac{n-1}{2}-|s|$, the composition factor $V_i(s)$  exists. It has the following $K$-types
$$\{ (\lambda, \lambda) \mid \lambda_i \geq |s|-\frac{n-1}{2}+i \geq \lambda_{i+1} \}.$$
These composition factors are called small constituents.
\item 
The composition factors $V_{\frac{n}{2}}(s)^{\pm}$ always exist. They are called large constituents. They have the following $K$-types:
$$\{ (\lambda, \lambda) \mid \pm \lambda_{\frac{n}{2}} \geq |s|+\frac{1}{2} \}.$$
\end{enumerate}
\item $V_i(s)$ is unitary only in the following two situations,
\begin{enumerate}
\item If $\frac{n-1}{2}-|s| \in [0, [\frac{n-1}{2}]]$, $V_{\frac{n-1}{2}-|s|}(s)$ is always unitary. 
\item If $n$ is even, $V_{\frac{n}{2}}(s)^{\pm}$ is always unitary for any $s \in \frac{n-1}{2}+\mathbb Z$, i.e., $s$ a half integer. $V_{\frac{n}{2}}(s)^+ \oplus V_{\frac{n}{2}}(s)^-$ yields an irreducible unitary representation of $SO(n,n)$.
\end{enumerate}
\item If $m=\frac{n-1}{2}-s$ and $[\frac{n-1}{2}] \geq m \geq 0$, $V_{m}(s)$ occurs as the unique quotient of $I_n(s)$ and as the unique subrepresentation of $I_n(-s)$. 
\end{enumerate}
\end{Theorem}
\subsection{Langlands-Vogan Parameter and Growth of $V_m(\frac{n-1}{2}-m)$}
Now suppose $s = \frac{n-1}{2}-m$ and $m \in [0, [\frac{n-1}{2}]]$.
Clearly, Vogan's lowest $K$-type of $V_i(s) \, (i \geq m)$ is $(\lambda, \lambda)$ with
$$\lambda=(\lambda_1= \ldots =\lambda_i(=-m+i) \geq \lambda_{i+1} = \lambda_{i+2}= \ldots =\lambda_{[\frac{n}{2}]}(=0)).$$
In particular, when $i=m$, the lowest $K$-type is trivial. Hence the small constituent $V_{m}(\frac{n-1}{2}-m)$ is spherical. We obtain
\begin{Theorem}\label{vm_estimate} Suppose $s = \frac{n-1}{2}-m$ and $m \in [0, [\frac{n-1}{2}]]$. Then
$V_m=V_{m}(\frac{n-1}{2}-m)$ has the following properties:
\begin{enumerate}
\item $V_m$ is spherical. The $K$-types of $V_m$ are of the form $(\lambda, \lambda)$ with
$$\lambda=(\lambda_1 \geq \lambda_2 \geq \ldots \geq \lambda_{m}(\geq 0)=\lambda_{m+1}= \ldots= \lambda_{\frac{n}{2}}).$$
In particular, if $m=0$, we obtain the \lq\lq smallest \!\rq\rq small constituent, the trivial representation.
\item The infinitesimal character $\mc I(V_m)=\eta(n-m-1, m)$ where
$$\eta(n-m-1, m)=(\overbrace{n-m-1,n-m-2, \ldots, m}^{n-2m},m, m-1,m-1, \ldots, 1,1, 0).$$
\item $V_m$ is the Langlands quotient  of $\Ind_{P_{\rm m}}^{G}  {\rm triv} \otimes \mathbb C_{ \eta(n-m-1,m)}$.
\item $V_m$ is the Vogan subquotient of $\Ind_{GL(n)_0 N}^{SO_0(n,n)} {\rm triv} \otimes  \|\det\|^{\frac{n-1}{2}-m}$.
\item The $K$-finite matrix coefficients of $V_m$ and of $I_n(\frac{n-1}{2}-m)$ are all bounded by a multiple of $\Xi_{\eta(n-m-1, m)}(g)$, and $\Xi_{\eta(n-m-1, m)}(\exp H^+)$ is bounded by
 $$(1+ (H^+, H^+))^Q \exp (\overbrace{-m,-m, \ldots, -m}^{n-2m},\overbrace{-m+1, -m+1}^2,-m+2,-m+2 \ldots, -1,-1,0,0)(H^+)$$
 for some positive $Q$.
\end{enumerate}
\end{Theorem}
\begin{cor}\label{bounds4in}
Let $u, v \in I_n(\frac{n-1}{2}-m) (m \in [0, [\frac{n-1}{2}]])$ such that $u|_K$ and $v|_K$ are continuous. Then
the matrix coefficient $(I_n(\frac{n-1}{2}-m)(g) u, v)$ is bounded by  
$$\sup( |u(k)|, k \in K) \sup( |v(k)|, k \in K)\Xi_{\eta(n-m-1, m)}(g).$$
 \end{cor}
 Proof: Since  $\eta(n-m-1, m)$ is real, $\Xi_{\eta(n-m-1, m)}(g)$ is a positive function. On the compact picture, the function $|u(k)|$ and $|v(k)|$ are bounded by  multiples of the spherical vector--the constant function ${\bf 1}_K$. Since $\frac{n-1}{2}-m$ is real,
 \begin{equation}
 \begin{split}
 |(I_n(\frac{n-1}{2}-m)(g) u, v)| \leq & \sup( |u(k)|, k \in K) \sup( |v(k)|, k \in K) (I_n(\frac{n-1}{2}-m)(g) {\bf 1}_K, {\bf 1}_K) \\
 = & \sup( |u(k)|, k \in K) \sup( |v(k)|, k \in K) \Xi_{\eta(n-m-1, m)}(g).
 \end{split}
 \end{equation}
 $\Box$.
 
\subsection{Howe's Duality Correspondence and Small constituent $V_m(\frac{n-1}{2}-m)$}

The representation $I_n(s)$ of $SO_0(n, n)$ can be extended to a representation of $SO(n,n)$ in two ways. By abusing notation, we keep $I_{n}(s)$ to denote the representation $\Ind_{GL(n)N}^{SO(n,n)} |\det|^s$.
All $V_i(s)$ now inherits a $SO(n,n)$ action from $I_n(s)$, except that
$V_{\frac{n}{2}}^+(s) \oplus V_{\frac{n}{2}}^-(s)$ constitute an irreducible representation of $SO(n,n)$ when $n$ is even. So from now on, the small constituent $V_{m}(\frac{n-1}{2}-m)$ will carry a $SO(n,n)$ action. \\
\\
In fact, $V_{m}(\frac{n-1}{2}-m)$ can be made into a representation of $O(n,n)$ as follows. We define an action of $O(n) \times O(n)$ on each $SO(n) \times SO(n)$-type $(\xi) \otimes (\xi)$ via
$(\xi,+) \otimes (\xi, +)$. Then $V_{m}(\frac{n-1}{2}-m)$ becomes a representation of $O(n,n)$. Again, by abusing notation, use $V_{m}(\frac{n-1}{2}-m)$ to denote this  representation of $O(n,n)$.\\
\\
Howe's duality correspondence is a one-to-one correspondence between a subset of the admissible dual of $O(p,q)$ and a subset of the admissible dual of $\widetilde{Sp}_{2k}(\mb R)$ (\cite{howe}). Let $(O(p,q), Sp_{2k}(\mb R))$ be a dual reductive pair in $Sp_{2k(p+q)}(\mb R)$. Let $\widetilde{Sp}_{2k(p+q)}(\mb R)$ be the unique double covering of $Sp_{2k(p+q)}(\mb R)$. Let $\omega$ be the oscillator representation of $\widetilde{Sp}_{2k(p+q)}(\mb R)$. Let $\pi_1$ and $\pi_2$ be irreducible admissible representations of $O(p,q)$ and $\widetilde{Sp}_{2k}(\mb R)$. We say that
$\pi_1$ corresponds to $\pi_2$ if $\pi_1^{\infty} \otimes \pi_2^{\infty}$ occurs as a quotient of $\omega^{\infty}$ (\cite{howe}). This is Howe's correspondence, also called local theta correspondence. We denote it by $\theta$. To be more specific, we use $\theta(p,q;2k)$ to denote the correspondence from $O(p,q)$ representations to $\widetilde{Sp}_{2k}(\mb R)$ representations and $\theta(2k; p,q)$ to denote the correspondence from $\widetilde{Sp}_{2k}(\mb R)$ representations to $O(p,q)$ representations. Sometimes, we use $\omega(p,q;2k)$ to denote $\omega$ restricted to the coverings of $(O(p,q), Sp_{2k}(\mb R))$. \\
\\
By a theorem of Lee (\cite{ll}), the small constituent $V_m(\frac{n-1}{2}-m)$  corresponds to a one dimensional character under Howe's duality correspondence. By the theory of stable range of local theta correspondence $\theta(2m; n, n)({\rm triv}) \cong \omega(n,n;2m)^{\infty} \otimes_{Sp_{2m}(\mb R)} {\rm triv}$ (\cite{li}). Notice here we must have $m \leq \frac{n-1}{2}$. So the dual pair 
$(O(n,n), Sp_{2m}(\mb R))$ is in the stable range.

\begin{Theorem}[\cite{ll}, \cite{li}]\label{leehe} Suppose that the integer $0 \leq m \leq \frac{n-1}{2}$.
The representation $\theta(2m; n, n)(\rm triv)$ is infinitesimally equivalent to $V_m(\frac{n-1}{2}-m)$.
In addition, $$V_m(\frac{n-1}{2}-m) \cong \omega^{\infty}(n,n; 2m) \otimes_{Sp_{2m}(\mb R)} {\rm triv}.$$
\end{Theorem}
Notice that when $p+q$ is even, $\theta(p,q;2k)$ is a correspondence between representations of the linear group $O(p,q)$ and representations of the linear group $Sp_{2k}(\mb R)$. 

\subsection{Matrix Coefficients of the Oscillator Representation and $\Xi_{\eta(n-m-1,m)}(g)$}

Consider the oscillator representation $(\omega, L^2(\mathbb R^{N}))$ of $\widetilde{Sp}_{2N}(\mb R)$ as in \cite{he00} \cite{basic}. Choose the standard maximal compact subgroup $K$. Then $K$ is the double covering of $U(N)$. There is a cannonical unitary character $\sqrt{\det}$ on $K$. The Harish-Chandra module is  $\mc P(x) \exp -\frac{1}{2} \| x \|^2$ with $\mc P(x)$ the polynomial algebra on $x \in \mathbb R^N$. 
Let $\mc P_k(x)$ be homogeneous polynomials of degree $k$. Let $\mc P_{\leq k}(x)$ be the polynomials of degree less or equal to $k$. Then $\mc P_k(x) \exp -\frac{1}{2} \| x \|^2$ is not a $K$-invariant subspace, except when $k=0$. But $\mc P_{\leq k}(x) \exp -\frac{1}{2} \| x \|^2$ is a $K$-invariant subspace. $\omega$ is a lowest weight module with the lowest weight vector $\exp -\frac{1}{2} \| x \|^2$. The group $K$ acts on $\exp -\frac{1}{2} \| x\|^2$ by scalar $\sqrt{\det}$. We sometimes call
$\exp -\frac{1}{2} \| x\|^2$ an almost spherical vector.
\begin{lem}\label{basices} Write $g=k_1 \exp H(g) k_2$
with $k_1, k_2 \in K$ and $H=(H_1, H_2, \ldots H_N) \in \mathbb R^N$. Then
$$F_N(g)=(\omega(k_1 \exp H k_2) \exp -\frac{1}{2} \| x \|^2, \exp -\frac{1}{2} \| x \|^2)= C_N (\sqrt{\det}(k_1 k_2)) \prod_{i=1}^N (\exp H_i+ \exp -H_i)^{-\frac{1}{2}},$$
where $C_N$ is a positive constant depending only on $N$. Obviously, if $\sqrt{\det}(k_1 k_2) =1$, then $F_N(g) >0$. 
\end{lem}
Now let $N=2nm$. Embed the dual pair $(O(n,n), Sp_{2m}(\mb R))$ as in \cite{basic}. Essentially, this amounts to a compactible
$KAK$ decomposition for the dual pair in $Sp_{2N}(\mb R)$. In this particular situation, we actually have $(O(n,n), Sp_{2m}(\mb R))$ acting on $L^2(\mb R^{2n} \otimes \mb R^m)$. Let $k \in \widetilde U(m)$. Then $k$ acts on $\exp -\frac{1}{2} \| x\|^2$ by $\sqrt{\det}(k)^n \overline{\sqrt{\det}(k)}^{n}=1$.
Obviously the compact group $O(n)O(n)$ acts on $\exp -\frac{1}{2} \| x\|^2$  trivially. Hence we obtain
$$F_N|_{O(n,n) Sp_{2m}(\mb R)} > 0.$$
Let $1$ be a unit vector in the trivial representation of $Sp_{2m}(\mb R)$. By Howe's theory (\cite{howe}), Theorem \ref{leehe} and \cite{he00}, the invariant tensor $\exp -\frac{1}{2} \| x \|^2
\otimes_{Sp_{2n}(\mb R)} 1$ is a spherical vector of $V_m(\frac{n-1}{2}-m)$. Hence we obtain
\begin{Theorem}\label{xios} Let $ (O(n,n), Sp_{2m}(\mb R))$ be a dual pair in $Sp_{4nm}(\mb R)$ as in \cite{basic}. Then the matrix coefficient $$(\omega(g_1, g_2) \exp -\frac{1}{2} \| x \|^2, \exp -\frac{1}{2} \| x \|^2) > 0 \qquad ((g_1, g_2) \in (O(n,n), Sp_{2m}(\mb R)))$$
and 
$$\Xi_{\eta(n-m-1, m)}(g_1)=C \int_{g_2 \in Sp_{2m}(\mb R)} (\omega(g_1, g_2) \exp -\frac{1}{2} \| x \|^2, \exp -\frac{1}{2} \| x \|^2) d g_2.$$
\end{Theorem}
I shall add one well-known lemma concerning the matrix coefficients of the oscillator representation.
\begin{lem}\label{kfinitees} Let $(\omega, L^2(\mb R^N))$ be the oscillator representation of $\widetilde{Sp}_{2N}(\mb R)$. Let $\phi, \psi \in \mc P(x) \exp -\frac{1}{2} \| x \|^2$ be a $K$-finite vector. Then there exists a positive constant $K_{\phi, \psi}$ such that
$$|(\omega(k_1 \exp H k_2) \phi, \psi)| \leq K_{\phi, \psi} |F_N(k_1 \exp H k_2)| =K_{\phi, \psi} C_N \prod_{i=1}^N (|H_i^{-1}|+ |H_i|)^{-\frac{1}{2}}.$$
\end{lem}
Proof: Suppose that all $K$-translations of $\phi$ and $\psi$ are uniformly bounded by $C (1+\| x \|^2)^M \exp -\frac{1}{2} \| x \|^2$. This is possible because $K$-translations of $\phi$ and $\psi$ are in a compact subset in a finite dimension subspace of $\mc P(x) \exp -\frac{1}{2} \| x \|^2$. Clearly $C (1+\| x \|^2)^M \exp -\frac{1}{2} \| x \|^2$ is bounded by $C_1 \exp -\frac{1}{4} \| x \|^2$. An easy computation shows that 
$$|(\omega(k_1 \exp H k_2) \phi, \psi)| \leq C_1^2 (\omega(\exp H) \exp -\frac{1}{4} \| x \|^2, \exp -\frac{1}{4} \| x \|^2)=K_{\phi, \psi} | F_N(k_1 \exp H k_2)|.$$
$\Box$

\section{Quantum Induction and Positivity of the Invariant Hermitian Form }
Let $\pi$ be an admissible representation of a semisimple Lie group $H$. Fix a maximal compact subgroup $K_H$ of $H$. {\bf We use $V(\pi)$ to denote the Harish-Chandra module of $\pi$.} In most part of this paper, the maximal compact subgroup is specified explicitly or implicitly and $V(\pi)$ is well-defined. \\
\\
 {\bf Let $n=p+q+d$}.
Equip the vector space $\mathbb R^{n, n}=\mathbb R^{2n}$ with the standard quadratic form of signature $(n, n)$ and the standard basis $\{e_i, f_j \mid i \in [1,n], j \in [1, n] \}$. Decompose $\mathbb R^{n, n}$ canonically as $\mathbb R^{p,q} \oplus \mathbb R^{q+d, p+d}$. This induces a cannonical embedding of $SO(p,q) SO(q+d, p+d)$ into $SO(n, n)$. The maximal compact subgroups are chosen to be the standard ones. \\
\\
{\bf Let $0 \leq m \leq \frac{n-1}{2}$. Denote the unitary representation $\theta(2m;n,n)(\rm triv)$ by $\mc E_m(n)$}. Then $\mc E_m(n)$ is  infinitesimally equivalent to $V_{m}(\frac{n-1}{2}-m)$. We also use  $\mc E_m(n)$ to denote the group action and Lie algebra action.  {\bf Let $(\pi, \mc H_{\pi})$ be an irreducible unitary representation of $SO(p,q)$.  Suppose that $V(\mc E_m(n)) \otimes_{G} V(\pi)$ is well-defined and nonzero}.  The nonvanishing of $V(\mc E_m(n)) \otimes_{SO(p,q)} V(\pi)$ will be established in Theorem \ref{quan_theta}. \\
\\
From Lemma \ref{tensorrepn}, we see that $V(\mc E_m(n)) \otimes_{SO(p,q)} V(\pi)$ is a $(\f g, K)$-module for $SO(q+d, p+d)$. Recall from Prop. \ref{invher} that $V(\mc E_m(n)) \otimes_SO(p,q) V(\pi)$ is equipped with a canonical $(\f {o}(q+d, p+d), O(q+d) O(p+d))$ invariant Hermitian form
$$(\phi_1 \otimes_{{SO(p,q)}} u_1, \phi_2 \otimes_{{SO(p,q)}} u_2)=\int_{{SO(p,q)}} (\mc E_m(n)(g) \phi_1, \phi_2)(\pi(g) u_1, u_2) d  g.$$
Fix a nonzero vector $u \in V(\pi)$. In this section, we shall show that this Hermitian form on $V(\mc E_m(n)) \otimes_{SO(p,q)} u$  is positive definite.  Therefore, the $(\f {o}(q+d, p+d), O(q+d) O(p+d))$-module  $V(\mc E_m(n)) \otimes_{SO(p,q)} u$ has an invariant pre-Hilbert structure. 
\begin{defn}\label{quantuminduction}
 We call the process:
$\mc Q(2m): V(\pi) \rightarrow V(\mc E_m(n)) \otimes_{SO(p,q)} V(\pi)$
quantum induction. Here $\mc Q(2m)(V(\pi)$ is a $(\f g, K)$-module for $SO(q+d, p+d)$.
\end{defn}

In \cite{quan},  $\mc Q(2m)$ is written more completely as $\mc Q(p,q;2m;q+d, p+d)$ for the orthogonal groups $O(*,*)$.  Quantum induction can be defined for all classical groups (\cite{quan1}). 

\subsection{A Positivity Theorem}

 We say a unitary representation of a unimodular group $G$ is {\it almost square integrable} if its matrix coefficients with respect to a dense subspace are in $L^{2+\epsilon}(G)$ for any $\epsilon > 0$.
\begin{Theorem}\label{posfn}
Let $G_1$ be a unimodular subgroup of a semisimple Lie group $G$. Let $(\sigma, \mc H_{\sigma})$ be an almost square integrable representation of $G$ and $(\pi, \mc H_{\pi})$ be an unitary representation of $G_1$. Here $\sigma$ is often not irreducible.  Let $K$ be a maximal compact subgroup of $G$.   Let $\Xi$ be Harish-Chandra's basic spherical function for $G$. Let $\psi \in \mc H_{\pi}$ such that
$$\int_{G_1} | \Xi(g_1)(\pi(g_1) \psi, \psi) | d g_1 < \infty.$$
 Then the invariant tensor product $(\mc H_{\sigma})_K \otimes_{G_1} \mathbb C \psi$ is well-defined and its canonical invariant Hermitian form is positive definite.
\end{Theorem}
Proof: Fix any $\phi \in (\mc H_{\sigma})_K$. Since $\sigma$ is almost square integrable, by a Theorem of Cowling-Haagerup-Howe (\cite{chh}), for any $\phi_1 \in (\mc H_{\sigma})_K$, there is a constant $C$ such that
$| (\sigma(g) \phi, \phi_1)| \leq C \Xi(g).$
Since $\int_{G_1} | \Xi(g_1)(\pi(g_1) \psi, \psi) | d g_1 < \infty$, 
$(\sigma(g_1) \phi, \phi_1)(\pi(g_1) \psi, \psi) \in L^1(G_1)$. So $(\mc H_{\sigma})_K \otimes_{G_1} \mathbb C \psi$ is well-defined. Now it suffices to show that $(\phi \otimes_{G_1} \psi, \phi \otimes_{G_1} \psi) \geq 0$. \\
\\
 Since $\phi$ is $K$-finite, let $S \subset \hat K$ be the finite set of $K$-types occuring in  the $K$-module generated by $\phi$. Let $C_c^{\infty}(G)(S)$ be the subspace of $C_c^{\infty}(G)$ such that the $K$-action from the left decomposes into the $K$-types in $S$. By our assumption, $(\sigma(g) \phi, \phi)$ is almost square integrable.
So $(\sigma, \mc H_{\sigma})$ is weakly contained in $L^2(G)$ (\cite{chh}). We can construct a sequence of $K$-finite functions $\phi_i$ in $C_c^{\infty}(G)(S)$ such that $(\sigma(g) \phi, \phi)$ can be approximated by $(L(g) \phi_i, \phi_i)$ uniformly on compacta. 
In particular, taking $g$ be the identity, we obtain
$$\| \phi_i\|_{L^2} \rightarrow \| \phi \|_{\mc H_{\sigma}}.$$
Hence $\{ \| \phi_i \|_{L^2} \}$ is bounded uniformly for all $i$. 
Now by a Theorem of Cowling-Haagerup-Howe (\cite{chh}),
$$\|(L(g) \phi_i, \phi_i) \| \leq C \Xi(g) \qquad (\forall \,\, g \in G),$$
where $C$ only depends on the norm $\| \phi_i \|$ and the $K$-types in $S$. Thus $C$ can be chosen uniformly for all $ \phi_i$.
Observe that
\begin{equation}
\begin{split}
 & \int_{G_1} (L(g_1) \phi_i, \phi_i)(\pi(g_1) \psi, \psi) d g_1 \\
 = & \int_{G_1} \int_{G} \phi_i(g_1^{-1} g) \overline{\phi_i(g)} d g (\pi(g_1) \psi, \psi) d g_1 \\
 & [ absolutely \ \ integrable, \ \ since \ \ \phi_i \ \ is \ \ compactly \ \ supported \ ] \\
 = & \int_{G_1} \int_{G_1 \backslash G} \int_{G_1} \phi_i(g_1^{-1} h_1 g) \overline{\phi_i(h_1 g)} d h_1 d [ G_1 g] (\pi(g_1) \psi, \psi) d g_1 \\
 = & \int_{G_1 \backslash G} \int_{G_1 \times G_1} \phi_i(g_1^{-1} h_1 g) \overline{\phi_i(h_1 g)} (\pi(g_1) \psi, \psi) d g_1 d h_1 d [ G_1 g] \\
 = & \int_{G_1 \backslash G} \int_{G_1 \times G_1} \phi_i( g_1 g) \overline{\phi_i(h_1 g)} (\pi(h_1 g_1^{-1})\psi, \psi) d g_1 d h_1 d [G_1 g] \\
 = & \int_{G_1 \backslash G} \int_{G_1 \times G_1} \phi_i( g_1 g) \overline{\phi_i(h_1 g)} (\pi(g_1^{-1})\psi, \pi(h_1^{-1}) \psi) d g_1 d h_1 d [G_1 g] \\
 = & \int_{G_1 \backslash G} (\int_{G_1} \phi_i( g_1 g) \pi(g_1^{-1})\psi d g_1, \int_{G_1}  \phi_i(h_1 g) \pi(h_1^{-1}) \psi d h_1 ) d [G_1 g] \\
 = & \int_{G_1 \backslash G} \| \int_{G_1} \phi_i( g_1 g) \pi(g_1^{-1})\psi d g_1 \|^2_{\mc H_{\sigma}} d [G_1 g] \\
 \geq & 0
\end{split}
\end{equation}
Now $\{ |(L(g_1) \phi_i, \phi_i)(\pi(g_1) \psi, \psi)| \}$ are uniformly bounded by an integrable function $| C \Xi (g_1) (\pi(g_1) \psi, \psi) |$. By the Dominated Convergence Theorem, 
\begin{equation}
\begin{split}
\int (\sigma(g_1) \phi, \phi) (\pi(g_1) \psi, \psi) d g_1 = & \int \lim_{i \rightarrow \infty} (L(g_1) \phi_i, \phi_i)(\pi(g_1) \psi, \psi) d g_1 \\
= & \lim_{i \rightarrow \infty} \int_{G_1} (L(g_1) \phi_i, \phi_i)(\pi(g_1) \psi, \psi) d g_1 \geq 0.
\end{split}
\end{equation}
Hence $(\phi \otimes_{G_1} \psi, \phi \otimes_{G_1} \psi) \geq 0$. $\Box$

\subsection{Positivity of the Hermitian Form}
Now we shall apply Theorem \ref{posfn} and prove that under a certain growth condition on $\pi \in \Pi_u(SO(p,q))$ the canonical Hermitian form on  $V(\mc E_m(n)) \otimes_{SO(p,q)} u$ is positive definite.  \\
\\
By Theorem \ref{vm_estimate}, the matrix coefficients of $V(\mc E_m(n))$ are bounded by a multiple of
$$(1+ (H^+(g), H^+(g)))^Q \exp (\overbrace{-m,-m, \ldots, -m}^{n-2m},\overbrace{-m+1, -m+1}^2,\overbrace{-m+2,-m+2}^2 \ldots, -1,-1,0,0)(H^+(g))$$
Since $n-2m \geq 1$, their restrictions onto $SO(k, 2m+2-k)$ are almost square integrable. Hence $\mc E_m(n)|_{SO(k, 2m+2-k)}$ is almost square integrable. We can now apply Theorem \ref{posfn} with  $G=SO(k, 2m+2-k)$ and $G_1=SO(p,q)$. The necessary and sufficient condition for the existence of a $G_1 \subseteq G$ is that $p+q \leq 2m+2$.
\begin{Proposition}\label{posinv}
Suppose that $p+q \leq 2m+2 \leq n+1$ and $p \leq q$. Let $\pi$ be an irreducible unitary representation of $SO(p,q)$ such that
its every $K$ finite matrix coefficient $f(g)$ satisfies the condition that
$$ |f(g)| \leq C_f \exp (\overbrace{m+2-p-q- \epsilon, m+3-p-q, \ldots m+1-q}^{p})(| H^+(g)|).$$
for some $\epsilon >0$. Then for any $u \in V(\pi)$, the invariant tensor product $V(\mc E_{m}(n)) \otimes_{SO(p,q)} u$ is well-defined and the canonical Hermitian form on $V(\mc E_{m}(n)) \otimes_{SO(p,q)} u$ is positive definite.
\end{Proposition}
Proof: Suppose that $f(g)=(\pi(g)u, u)$ satisfies the condition that
$$ |f(g)| \leq C \exp (\overbrace{m+2-p-q- \epsilon, m+3-p-q, m+4-p-q, \ldots, m+1-q}^{p})(| H^+(g)|).$$
Recall that $\Xi(x)$ for the group $SO(p, 2m+2-p)$ is bounded by
$$ C_1 (1+(H^+(x), H^+(x)))^Q \exp(\overbrace{-m, -m+1, -m+2, \ldots, -m+p-1}^{p})(H^+(x)).$$
Then $\Xi|_{SO(p,q)}(g) |f(g)|$ is bounded by
$$C_1 C (1+(H^+(g), H^+(g)))^Q \exp(\overbrace{2-p-q-\epsilon, 4-p-q, 6-p-q, \ldots, 0}^{p})(|H^+(g)|),$$
which is clearly in $L^1(SO(p,q))$. So $\Xi|_{SO(p,q)}(g) |f(g)|$ is in $L^1(SO(p,q))$.  We have seen that $\mc E_m(n)|_{SO(p, 2m+2-p)}$ is almost square integrable.
By Theorem \ref{posfn},  the invariant tensor product $V(\mc E_{m}(n)) \otimes_{SO(p,q)} u$ is well-defined. In addition, for any $\phi \in V(\mc E_m(n))$, we have
$$(\phi \otimes_{SO(p,q)} u, \phi \otimes_{SO(p,q)} u) \geq 0.$$
Hence the invariant Hermitian form on $V(\mc E_{m}(n)) \otimes_{SO(p,q)} u$ is positive definite.
$\Box$  \\
\\
Under the hypothesis of Prop. \ref{posinv}, there is a natural injection $i$ from $V(\mc E_{m}(n)) \otimes_{SO(p,q)} \mathbb C u$ to $V(\mc E_{m}(n)) \otimes_{SO(p,q)} V(\pi)$.  This can be seen as follows. By Cor \ref{bounds4in} and Theorem \ref{l1dominatedTheorem}, we have
$$V(\mc E_{m}(n)) \otimes_{SO(p,q)} \mathbb C u \cong V(\mc E_m(n)) \otimes_{(\mc E_m(n)^{\infty})^c \otimes \mathbb C u} \mathbb C u;$$
$$V(\mc E_{m}(n)) \otimes_{SO(p,q)} V(\pi) \cong V(\mc E_{m}(n)) \otimes_{(\mc E_m(n)^{\infty})^c \otimes V(\pi)^c} V(\pi).$$
Obviously, if $\phi \otimes_{SO(p,q)} u=0$ in $V(\mc E_{m}(n)) \otimes_{\mc E_m(n)^c \otimes V(\pi)^c} V(\pi)$, then $\phi \otimes_{(\mc E_m(n))^c \otimes \mathbb C u} u=0$.
Suppose that $\phi \otimes_{SO(p,q)} u=0$ in $V(\mc E_{m}(n)) \otimes_{SO(p,q)} u$. This means that for any $\phi_1 \in \mc E_{m}(n)^{\infty}$, we have
$$\int_{SO(p,q)} (\mc E_m(n) (g)) \phi, \phi_1 )(\pi(g) u,  u) d g=0.$$
It follows that for any $\phi_1 \in \mc E_{m}(n)^{\infty}, g_1 \in SO(p,q)$,
$$0=\int_{SO(p,q)} (\mc E_m(n)(g_1^{-1} g) \phi,  \phi_1) (\pi(g_1^{-1}g) u, u) d g=\int_{SO(p,q)} (\mc E_m(n)( g) \phi, \mc E_m(n)(g_1) \phi_1) (\pi(g) u, \pi(g_1)u) d g.$$
Let $W$ be the finite linear span of  $\{ \pi(g_1) u \mid g_1 \in SO(p,q) \}$.
Since $\mc E_m(n)(g_1) \phi_1 \in \mc E_m(n)^{\infty}$, we have $\phi \otimes_{(\mc E_m(n)^{\infty})^c \otimes W^c} u=0$. Notice that $(\mc E_m(n)^{\infty}) \otimes W$ and $(\mc E_m(n)^{\infty}) \otimes V(\pi)$ are mutually $L^1$ bounded. Again by Theorem \ref{l1dominatedTheorem}, $\phi \otimes_{(\mc E_m(n)^{\infty})^c \otimes V(\pi)^c} u=0$. So $\phi \otimes_{SO(p,q)} u=0$ in $V(\mc E_m(n)) \otimes_{SO(p,q)} V(\pi)$. Thus the natural map $$i: V(\mc E_{m}(n)) \otimes_{SO(p,q)} \mathbb C u \rightarrow V(\mc E_{m}(n)) \otimes_{SO(p,q)} V(\pi)$$
is injective. 

\subsection{Independence of the Hilbert Structure}
Let $s=\frac{n-1}{2}-m \geq 0$.
Recall that $V_m(-s)$ is a subrepresentation of $I_n(-s)$ and it inherits a pre-Hilbert space structure and $(\f o(n,n), S(O(n)O(n)))$-module structure from $I_n(-s)$. By Lemma \ref{gkh} and Remark \ref{gkh1}, $V_m(-s) \otimes_{SO(p,q)} V(\pi)$ inherits a $(\f o(q+d,p+d), S(O(q+d)O(p+d)))$-module structure from $V_{m}(-s)$. 
We have seen that on the Harish-Chandra module level $\mc E_m(n)$ can be identified with $V_m(-s)$. But their inner products are different. So given a Hilbert representation $(\pi, \mc H)$ of $G$, can we identifty $V(\mc E_m(n)) \otimes_{SO(p,q)} V(\pi)$ with
$V_m(-s) \otimes_{SO(p,q)} V(\pi)$ on the Harish-Chandra module level? \\
\\
Let $j: V(\mc E_m(n)) \rightarrow V_m(-s)$ be the identification of the Harish-Chandra module. Notice that both inner products on $V(\mc E_m(n))$ and $V_m(-s)$ are $K$-invariant.   Let $( \, , \,)_1$ be the inner product on $V(\mc E_m(n))$ and $(\, , \,)_2$ be the inner product on $V_m(-s)$. We can construct a map $A: V(\mc E_m(n)) \rightarrow V_m(-s)$ such that for any $\phi, \psi \in V(\mc E_m(n))$
$$(\phi, \psi)_1=(j(\phi), A (\psi))_2.$$
Since both inner products are nondegenerate, $A$ must be a bijection on each $K$-type. So $A$ is one-to-one and onto. We have
$$(\mc E_m(n)(g) \phi,  \psi)_1=(I_n(-s)(g) j(\phi), A (\psi))_2.$$
In particular $\sum_{i} \phi_i \otimes_{SO(p,q)} v_i \in V(\mc E_m(n)) \otimes_{SO(p,q)} V(\pi)$ vanishes if and only if 
$$\int_{SO(p,q)}((\mc E_m(n)(g) \phi_i,  \psi)_1(\pi_i(g) v_i, u) d g =0 \qquad (\forall \, \, \psi \in V(\mc E_m(n)), u \in V(\pi)),$$ 
if and only if
$$\int_{SO(p,q)}((I_n(-s)(g) j(\phi_i), A(\psi))_2(\pi_i(g) v_i, u) d g =0 \qquad (\forall \, \, \psi \in V(\mc E_m(n)), u \in V(\pi)).$$ 
We see that the kernels of the two averaging operators are the same. We obtain

\begin{Theorem}\label{emvm} Let  $n=p+q+d$ and $2m+1 \leq n$. Let $(\pi, \mc H)$ be an irreducible admissible Hilbert representation of $G=SO(p,q)$. As Harish-Chandra modules of $SO(q+d, p+d)$, we have $V(\mc E_m(n)) \otimes_G V(\pi) \cong V_m(m-\frac{n-1}{2}) \otimes_G V(\pi)$. Suppose that there is an {\it invariant} Hermitian form attached to the smooth vectors $\mc H^{\infty}$. Let ${}^h V(\pi)$ be the Harish-Chandra module equipped with this Hermitian structure. Then as Harish-Chandra modules of $SO(q+d, p+d)$,
$$V(\mc E_m(n)) \otimes_G V(\pi) \cong V(\mc E_m(n)) \otimes_G {}^h V(\pi).$$
\end{Theorem}
The second assertion follows from the same argument as the first assertion.
I shall also remark that this theorem holds for any semisimple group $G$ and any irreducible unitary Harish-Chandra module $V(\mc E)$ of a bigger semisimple group $Q \supseteq G$.
 
\section{Quantum Induction: Subrepresentation Theorem}
{\bf Let  $n=p+q+d$, $0 \leq m \leq \frac{n-1}{2}$ and  $s=\frac{n-1}{2}-m$. Let $\eta=\eta(n-m-1, m)$.}
In this section, we shall  show that $\mc Q(2m)(V(\pi))$ is a subrepresentation of the induced module $\Ind_{SO(q,p) GL(d) N}^{SO(q+d, p+d)} \pi \otimes |\det|^{-\frac{n-1}{2}+m}$ after we identify $SO(q, p)$ with $SO(p, q)$. Hence $\mc Q(2m)(V(\pi))$ is admissible and quasisimple.
\\
\\
Let us consider the $SO(p,q) SO(q+d, p+d)$ action on $I_n(s)$.
$I_{n}(s)$ is a quasisimple admissible representation of $SO(n, n)$. The $K$-finite subspace  $V(I_{n}(s))$ is a Harish-Chandra module. Vectors in $I_{n}(s)$ can be regarded as functions on $X \cong SO(n,n)/P_n$, one branch of the maximal isotropic Grassmanian of $\mathbb R^{n,n}$.  Here $P_n$ is the Siegel parabolic subgroup.
The symmetric subgroup $SO(p,q) SO(q+d, p+d)$ acts on $X$ with a unique open dense orbit $X_0$ (\cite{hh}). The action on this open dense orbit can be described as follows. \\
\\
Given two $n$-dimensional Euclidean spaces $\mathbb R^n$ equipped with the standard inner products,
let $\mathbb R^{n,n}= \mathbb R^n \oplus \mathbb R^n$ be a $2n$-dimensional real vector space equipped with the symmetric form 
$$((x_1, y_1), (x_2, y_2))=(x_1, x_2)-(y_1, y_2), \qquad (x_1, x_2, y_1, y_2 \in \mathbb R^n).$$
Fix the standard basis $\{ e_1, e_2, \ldots, e_n, f_1, f_2, \ldots, f_n \}$. Let $SO(n,n)$ be the special orthogonal group preserving the symmetric form $(\, , \,)$. We decompose
$\mathbb R^{n,n}= \mathbb R^{p,q} \oplus \mathbb R^{q+d, p+d} $ where
$$\mathbb R^{p,q} = {\rm span}(e_1, e_2, \ldots, e_p, f_1, f_2 \ldots, f_q), \qquad 
\mathbb R^{q+d, p+d}={\rm span}(e_{p+1}, e_{p+2}, \ldots e_n, f_{q+1}, f_{q+2}, \ldots, f_{n}). $$ 
Let $SO(p,q)$ be the special orthogonal group preserving $(\, , \,)|_{\mathbb R^{p,q}}$, and $SO(q+d, p+d)$ be the special orthogonal group preserving
$(\, , \,)|_{\mathbb R^{q+d, p+d}}$.
We obtain a diagonal  embedding from $SO(p,q) SO(q+d, p+d)$ into $SO(n,n)$. Let
 $$W_1={\rm span}\{ e_1+f_{q+1}, e_2+f_{q+2}, \ldots, e_p+f_{q+p}; f_1+e_{p+1}, f_2+e_{p+2}, \ldots, f_q+e_{p+q}\}$$
 $$W_0={\rm span} \{ e_{p+q+1}+f_{p+q+1}, \ldots, e_{p+q+d}+f_{p+q+d} \}.$$
 \begin{lem}[\cite{hh}]\label{basepoint}
 \begin{enumerate}
 \item
Let $W=W_1 \oplus W_0$. Then $W$ is a maximal isotropic subspace of $\mathbb R^{n,n}$. 
\item Fix $W$ as the base point in  $ X$.
 Let $LN$ be the parabolic subgroup of $SO(q+d, p+d)$ that preserves $W_0$.
 Then $L \cong SO(q,p) GL(d)$ and $X_0 \cong SO(q+d, p+d)/GL(d) N$.
 \item $W$ induces a negative isometry between
 $\mathbb R^{p,q}$ and $\mathbb R^{q, p} \subseteq \mathbb R^{q+d, p+d}$:
  $$e_i \leftrightarrow f_{q+i}, \, \,  (i \in [1,p]); \qquad  f_j \leftrightarrow e_{p+j}, \,\,
  (j \in [1, q]).$$
 \item {\bf  $W$  induces an identification of $SO(p,q)$ with $SO(q, p)$, denote this by
 $g_1 \rightarrow \dot{g_1}$}.  
 \item Then $(SO(p,q) SO(q+d, p+d))_W= \{(g_1, \dot{g_1}) \mid g_1 \in SO(p,q) \} GL(d) N.$
 \end{enumerate}
 \end{lem}
{\bf From now on, $SO(p,q)$ will be identified with $SO(q,p)$ via $g \rightarrow \dot{g}$.} Let $|\det|$ be the character of $P_n$ obtained by taking the absolute value of the determinant character on the Levi factor. The restriction of $|\det|$ onto  $(SO(p,q) SO(q+d, p+d))_W$ is simply the absolute value of the determinant character of the $GL(d)$ factor. \\
 \\
 {\bf Let $C_c^{\infty}(X_0, s)$ be smooth sections in $I_{n}(s)$ with compact support in $X_0$.} Notice that
 both $C_c^{\infty}(X_0, s)$ and $V(I_n(s))$ are dense subspaces of the Hilbert representation $I_n(s)$.
 \subsection{Realization of $C_c^{\infty}(X_0, s) \otimes_{SO(p,q), C_c^{\infty}(X_0, -s) \otimes V(\pi^*)} V(\pi)$} The main result of this subsection is Theorem \ref{itpcompact}. This is the vector bundle version of Theorem \ref{georeal}. \\
 \\
  Recall that there is a $SO(n,n)$-invariant complex linear pairing between $I_n(s)$ and $I_n(-s)$:
  \begin{equation}\label{integraloverk}
 (\phi, \psi)=\int_{S(O(n) O(n))/O(n)} \phi(k) \psi(k) d [k], \qquad (\phi \in I_n(s), \psi \in I_n(-s)).
 \end{equation}
 This pairing induces a pairing between $C_c^{\infty}(X_0, s)$ and $C_c^{\infty}(X_0, -s)$.
 $C_c^{\infty}(X_0, s)$, regarded as a smooth representation of $SO(p,q) SO(q+d, p+d)$, consists of smooth sections of the homogeneous line bundle:
 \begin{equation}\label{bundle}
 SO(q+d, p+d) \times_{GL(d) N} |\det|^{\frac{n-1}{2}+s} \rightarrow X_0,
 \end{equation}
 with compact support.  \\
 \\
 Let $g_1 \in SO(p,q)$ and $x \in X_0$. By Lemma \ref{basepoint}, $x$ can be written as $g_2 W$ with $g_2 \in SO(q+d, p+d)$. Then 
 $$ g_1 g_2 W=g_2 g_1 W=g_2 g_1 g_1^{-1} \dot{g_1}^{-1} W= g_2 \dot{g_1}^{-1} W.$$
 Hence left $g_1$ action on $X_0 \cong SO(q+d, p+d) /{GL(d) N}$ is simply the right $\dot{g_1}^{-1}$ action. In particular
 $$(I_n(s)(g_1) f)(g_2)= f(g_1^{-1} g_2)= f(g_2 \dot{g_1})=(R(\dot{g_1}) f )(g_2).$$
  Now {\bf  equip $X_0$ with a left $S(O(q+d) O(p+d))$ and right $SO(q,p)$ invariant measure}. From the homogeneous line bundle structure $(\ref{bundle})$, we see that 
  the action of $S(O(q+d)O(p+d)) \times SO(q,p)$ on $X_0$ is transitive.  The stabilizer $T$ is isomorphic to $O(d) \Delta(S(O(q) O(p)))$. Here
  $$\Delta(S(O(q) O(p)))= \{ (k, k) \mid k \in S(O(q) O(p)) \} \subseteq S(O(q+d) O(p+d)) \times SO(q, p).$$
  Applying coordinate transformation to Eq. (\ref{integraloverk}), we have
  $$(\phi, \psi)= \int_{S(O(q+d)O(p+d)) \times SO(q,p)/ T} \phi(g) \psi(g) J([g]) d [g], \qquad (\phi \in I_n(s), \psi \in I_n(-s)) .$$
  Here $J([g])$ is the Jacobian of the coordinate transformation from
  $$ X_0 \ni [k]  \rightarrow [g] \in S(O(q+d)O(p+d)) \times SO(q,p)/ T.$$
 Notice that $\phi(g)\psi(g)$ only depends on $[g]$. Since this form is invariant under the action of $S(O(q+d)O(p+d)) \times SO(q,p)$, $J([g])$ must be a constant. We normalize the measure on $X_0$ so that $J([g])=1$. \\
 \\ 
  {\bf We identify $SO(p,q)$ with $SO(q,p) \subseteq SO(q+d, p+d)$ by $g \rightarrow \dot{g}$}. Let $(\pi, \mc H_{\pi})$ be an irreducible unitary representation of $SO(p,q)$.
 By Lemma \ref{pairing}, for any $u \in V(\pi)$, $v \in V(\pi^*)$, $\phi \in C_c^{\infty}(X_o, s)$ and $\psi \in C_c^{\infty}(X_0, -s)$,  we have
 \begin{equation}
 \begin{split}
  & (\phi \otimes_{SO(p,q), C_c(X_0, -s) \otimes V(\pi^*)} u, \psi \otimes_{SO(p,q), C_c(X_0, s) \otimes V(\pi)} v) \\ 
  = & \int_{SO(p,q)} (I_n(s)(g) \phi, \psi)(\pi(g) u, v) d g \\
  = & \int_{SO(q,p)} \int_{S(O(q+d)O(p+d)) \times SO(q,p)/T} \phi(x g) \psi(x) d [x] (\pi(g) u, v) d g \\
  = & \int_{SO(q,p)} \int_{(S(O(q+d)O(p+d)) \times SO(q,p))/ (T  SO(q,p))} \int_{SO(q,p)}\phi(x h g) \psi(x h) d h \, d [x] (\pi(g) u, v)  d g \\
  = & \int_{(S(O(q+d)O(p+d)) \times SO(q,p))/ (T  SO(q,p))} \int_{SO(q,p) \times  SO(q, p)} \phi(x h g)  \psi(x  h) (\pi(g) u, v) d g d h d [x] \\
  = & \int_{S(O(q+d) O(p+d))/ O(d) S(O(q)O(p))} \int_{SO(q,p) \times SO(q,p)} \phi(x h^{\prime}) \psi(x h)
  (\pi(h^{-1}) \pi(h^{\prime}) u, v) d h d h^{\prime} d [x] \\
  = & \int_{S(O(q+d) O(p+d))/ O(d) S(O(q)O(p))} (\int_{SO(q,p)} \phi(x h^{\prime}) \pi(h^{\prime}) u d h^{\prime}, \int_{SO(q,p)} \psi(x h) \pi^*(h) v d h) d [x]
 \end{split}
 \end{equation}
 It is easy to check that all these integrals are well-defined and converge absolutely. {\bf Now put
 $$\mc I(\phi \otimes u)(x)=\int_{SO(q,p)} \phi(x h) \pi(h) u d h \qquad (x \in SO(q+d, p+d)).$$}
 The map
 $\mc I$ is essentially integration on the right $SO(q,p)$ fiber for
 $$SO(q,p) \rightarrow SO(q+d, p+d)/GL(d) N \rightarrow SO(q+d, p+d)/SO(q,p) GL(d)N.$$
 From the structure of the bundle (\ref{bundle}), 
 $\mc I(\phi \otimes u))$
  can be identified with a smooth section of
 $$ SO(q+d, p+d) \times_{SO(q,p) GL(d) N} (\mc H_{\pi} \otimes |\det|^{\frac{n-1}{2}+s}) \rightarrow SO(q+d, p+d)/SO(q, p)GL(d) N.$$ 
 Similarly,
 $\mc I(\psi \otimes v))$
  can be identified with a smooth section of
 $$ SO(q+d, p+d) \times_{SO(q,p) GL(d) N} (\mc H_{\pi^*} \otimes |\det|^{\frac{n-1}{2}-s}) \rightarrow SO(q+d, p+d)/SO(q, p)GL(d) N.$$
 In addition, the complex linear pairing $(\phi \otimes_{SO(p,q), C_c(X_0, -s) \otimes V(\pi^*)} u, \psi \otimes_{SO(p,q), C_c(X_0, s) \otimes V(\pi)} v)$ is exactly
 $$(\mc I(\phi \otimes u), \mc I(\psi \otimes v))_{S(O(q+d) O(p+d))/ O(d) S(O(q)O(p))}.$$
 Notice here that $SO(q+d, p+d)/SO(q, p)GL(d) N \cong S(O(q+d) O(p+d))/ O(d) S(O(q)O(p))$. \\
 \\
 Recall that for $P=SO(q,p) GL(d) N \subseteq SO(q+d, p+d)$, the half sum of positive roots in $\Delta(\f o(q+d, p+d), \f a)$ is exactly $\frac{d(p+q+d-1)}{2}=\frac{d(n-1)}{2}$. This corresponds to the character $|\det|^{\frac{n-1}{2}}$ og $GL(d)$. Hence, we see that $\mc I(C_c^{\infty}(X_0, s) \otimes V(\pi))$ can be identified with certain smooth sections of
 $$\Ind_{SO(q,p) GL(d) N}^{SO(q+d, p+d)} \pi \otimes |\det|^s.$$
 
\begin{Theorem}\label{itpcompact} Regard $\pi$ also as a representation of $SO(q,p)$ by identifying $SO(q,p)$ with $SO(p,q)$ as in Lemma \ref{basepoint}. For any $u \in V(\pi)$, $v \in V(\pi^*)$, $\phi \in C_c^{\infty}(X_o, s)$ and $\psi \in C_c^{\infty}(X_0, -s)$, we have
$$(\phi \otimes_{SO(p,q), C_c(X_0, -s) \otimes V(\pi^*)} u, \psi \otimes_{SO(p,q), C_c(X_0, s) \otimes V(\pi)} v)=(\mc I(\phi \otimes u), \mc I(\psi \otimes v))$$
where the right hand side is the complex linear pairing between $\Ind_{SO(q,p) GL(d) N}^{SO(q+d, p+d)} \pi \otimes |\det|^s$ and $\Ind_{SO(q,p) GL(d) N}^{SO(q+d, p+d)} \pi^* \otimes |\det|^{-s}$.
In addition,
$\mc I$ induces an $SO(q+d, p+d)$-invariant isomorphism from
$$C_c^{\infty}(X_0, s) \otimes_{SO(p,q), C_c^{\infty}(X_0, -s) \otimes V(\pi^*)} V(\pi)$$ onto a dense subspace of smooth vectors
in $\Ind_{SO(q,p) GL(d) N}^{SO(q+d, p+d)} \pi \otimes |\det|^s$. 
\end{Theorem}
Proof: 
\begin{enumerate}
\item The first statement is clear from the argument proceeding the theorem.
\item Next, we show that 
$\{ \mc I(\phi \otimes u) \mid \phi \in C_c^{\infty}(X_0, s), u \in V(\pi) \}$
is dense in $\Ind_{SO(q,p) GL(d) N}^{SO(q+d, p+d)} \pi \otimes |\det|^s$. Recall that the vector bundle (\ref{bundle}) is a $GL(d)N$-principal bundle. Choose a local trivialization on an open $SO(q,p)$-invariant $\Omega \subset X_0$ for the vector bundle (\ref{bundle}). This amounts to choosing a local trivialization for the bundle
$$ SO(q+d, p+d) \times_{SO(q,p) GL(d) N} {\rm triv} \otimes |\det|^{\frac{n-1}{2}+s} \rightarrow 
SO(q+d, p+d)/SO(q,p) GL(d) N. $$
Then we have
$C_c^{\infty}(\Omega, s) \cong C_c^{\infty}(\Omega, \mathbb C)$. By Theorem \ref{georeal},
$\mc I(C_c^{\infty}(\Omega, \mathbb C) \otimes V(\pi))$ is dense in $L^2(\Omega \times_{SO(p,q)} \mc H_{\pi})$. By partition of unit, we obtain the desired result.  
\item Finally, we want to show that $\mc I$ induces an isomorphism from $C_c^{\infty}(X_0, s) \otimes_{SO(p,q), C_c^{\infty}(X_0, -s) \otimes V(\pi^*)} V(\pi)$ into $\Ind_{SO(q,p) GL(d) N}^{SO(q+d, p+d)} \pi \otimes |\det|^s$. Notice that
$\sum_{i=1}^t \phi_i \otimes_{SO(p,q), C_c^{\infty}(X_0, -s) \otimes V(\pi^*)} u_i=0$ if and only if
$$(\sum_{i=1}^t \phi_i \otimes_{SO(p,q), C_c^{\infty}(X_0, -s) \otimes V(\pi^*)} u_i)(\psi \otimes v)=0 \qquad (\forall \, \, \psi \in C_c^{\infty}(X_0, -s), v \in V(\pi^*));$$
if and only if
$(\mc I(\sum_{i=1}^t \phi \otimes u_i), \mc I(\psi \otimes v))=0$. By (2), $\mc I(\psi \otimes v)$ is dense in
$\Ind_{SO(q,p) GL(d) N}^{SO(q+d, p+d)} \pi^* \otimes |\det|^{-s}$. 
We have seen that $\sum_{i=1}^t \phi_i \otimes_{SO(p,q), C_c^{\infty}(X_0, -s) \otimes V(\pi^*)} u_i=0$ if and only if $\mc I(\sum_{i=1}^t \phi \otimes u_i)=0$. So $\mc I$ induces an isomorphism. It is easy to see that this isomorphism preserves the action of $SO(q+d, p+d)$.
\end{enumerate}
$\Box$ \\
\\
Given $\sigma \in \widehat{S(O(q+d) O(p+d))}$, let $P_{\sigma}$ be the projection of the Hilbert representation $\Ind_{SO(q,p) GL(d) N}^{SO(q+d, p+d)} \pi \otimes |\det|^s$ onto its $\sigma$-isotypic subspace, which will be finite dimensional. Let $C_c^{\infty}(X_0, s)_{\sigma}$ be the $\sigma$ isotypic subspace. Then 
$$P_{\sigma}\mc I(C_c^{\infty}(X_0, s) \otimes V(\pi))=\mc I(C_c^{\infty}(X_0, s)_{\sigma} \otimes V(\pi)).$$  
By Theorem \ref{itpcompact}, the left hand side must be the full finite dimensional $P_{\sigma}(\Ind_{SO(q,p) GL(d) N}^{SO(q+d, p+d)} \pi \otimes |\det|^s)$. Otherwise, $\mc I(C_c^{\infty}(X_0, s) \otimes V(\pi)))$ will not be dense in $\Ind_{SO(q,p) GL(d) N}^{SO(q+d, p+d)} \pi \otimes |\det|^s$. Hence
$$\mc I(C_c^{\infty}(X_0, s)_{\sigma} \otimes V(\pi))=P_{\sigma}(\Ind_{SO(q,p) GL(d) N}^{SO(q+d, p+d)} \pi \otimes |\det|^s).$$
{\bf Let $D(X_0, s)$ be the $S(O(q+d)O(p+d))$-finite subspace of $C_c^{\infty}(X_0, s)$.} We have
\begin{cor}\label{itpkfinite} $\mc I$ induces an equivalence of Harish-Chandra modules of $SO(q+d, p+d)$:
$$D(X_0, s) \otimes_{SO(q,p),C_c^{\infty}(X_0, -s) \otimes V(\pi^*) } V(\pi) \rightarrow 
V(\Ind_{SO(q,p) GL(d) N}^{SO(q+d, p+d)} \pi \otimes |\det|^s).$$
\end{cor}
 \subsection{Two Invariant Tensor Products}
 
\begin{Theorem}\label{equivalence} Let $s=\frac{n-1}{2}-m$, $n=p+q+d$ and $\eta=\eta(n-m-1, m)$. Let $\pi$ be an irreducible unitary representation of $SO(p,q)$ such that every leading exponent $\mu$ satisfies
$$\eta+(\mu, {\mathbf 0}_{q+d})-(d+1, d+2, \ldots, d+p, {\mathbf 0}_{q+d}) \prec 0.$$
Then $D(X_0, -s) \otimes V(\pi^*)$ and $V(I_n(s))^c \otimes V(\pi)^c$ are two mutually $L^1$-dominated subspaces of $I_n(-s) \otimes \mc H_{\pi}^*$, with respect to $V(I_n(s)) \otimes V(\pi)$. We have
\begin{enumerate}
\item as  $(\f o(q+d, p+d), S(O(q+d) O(p+d)))$-modules,
$$V(I_{n}(s)) \otimes_{SO(p,q)} V(\pi) \cong  V(I_n(s)) \otimes_{SO(p,q), D(X_0, -s) \otimes V(\pi^*)} V(\pi) .$$
\item as $SO(q+d, p+d)$ representations,
$$ D(X_0, -s) \otimes_{SO(p,q), C_c^{\infty}(X_0, s) \otimes V(\pi)} V(\pi)^* \cong D(X_0, -s) \otimes_{SO(p,q), V(I_n(s)) \otimes V(\pi)} V(\pi)^*.$$
\item There is a nondegenerate $(\f o(q+d, p+d), S(O(q+d) O(p+d)))$-invariant complex linear pairing $(\, , \,)$:
$$(V(I_n(s)) \otimes_{SO(p,q)} V(\pi), \, D(X_0, -s) \otimes_{SO(p,q), C_c^{\infty}(X_0, s) \otimes V(\pi)} V(\pi^*)) \longrightarrow \mathbb C.$$
\end{enumerate}
\end{Theorem}
Proof: Fix $K=S(O(n)O(n))$. For each $\phi \in D(X_0, -s)$, $\phi|_K$ is continuous and thus bounded on $K$. Let $\psi_i \in V(I_n(s))^c \subseteq I_n(-s)$ such that $\psi_i|_K$ approaches $\phi|_K$ uniformly. This is always possible by essentially Stone-Weierstrass Theorem. Then for any $\phi_0 \in V(I_n(s))$ or $\phi_0 \in D(X_0, s)$
$$(I_n(s)(g) \phi_0, \psi_i) \rightarrow (I_n(s)(g) \phi_0, \phi)$$
and by Cor. \ref{bounds4in}, all $|(I_n(s)(g) \phi_0, \psi_i)|$ are unformly bounded by 
$$\| \phi_0|_K \|_{\sup} (\sup_i \| \psi_i|_K \|_{\sup}) \Xi_{\eta}(g).$$
By our assumption on $\mu$ and Theorem \ref{leadingex}, 
$$\mu+ 2 \rho(SO(p,q))+(\eta- \rho(SO(n,n)))|_{\mathbb R^p \oplus \mathbf 0_{q+d}} \preceq 0.$$
So $|\Xi_{\eta}(g)|_{SO(p,q)}(\pi(g)u, v)|$ is in
$L^1(SO(p,q))$.
 By Defintion \ref{l1dominated} and Remark \ref{l1dominatedremark}, $D(X_0, -s) \otimes V(\pi^*)$ is $L^1$-dominated by $V(I_n(s))^c \otimes V(\pi)^c$ with respect to both 
$D(X_0, s) \otimes V(\pi)$ and $V(I_n(s)) \otimes V(\pi)$. Similarly, $V(I_n(s))^c \otimes V(\pi)^c$
is also $L^1$-dominated by $D(X_0, -s) \otimes V(\pi^*)$. \\
\\
By Theorem \ref{l1dominatedTheorem}, we have
 $$V(I_{n}(s)) \otimes_{SO(p,q)} V(\pi) \cong V(I_n(s)) \otimes_{SO(p,q), D(X_0, -s) \otimes V(\pi^*)} V(\pi).$$
 
 \commentout{Notice that the matrix coefficients with respect to $(V(I_n(s)), D(X_0, s)^c)$ are the same as the matrix coefficients with respect to
 $(V(I_n(s)), D(X_0, -s)$. Hence
 $$\ker(\mc A_{SO(p,q), D(X_0, s)^c \otimes V(\pi)^c})=\ker(\mc A_{SO(p,q), D(X_0, -s) \otimes V(\pi)^* }).$$
 We obtain
$$ V(I_n(s)) \otimes_{SO(p,q), D(X_0, s)^c \otimes V(\pi)^c} V(\pi) \cong V(I_n(s)) \otimes_{SO(p,q), D(X_0, -s) \otimes V(\pi)*} V(\pi) .$$}
By Lemma \ref{tensorrepn} and Lemma \ref{tensorrepn1}, both spaces inherits the same $(\f o(q+d, p+d), S(O(q+d)O(p+d)))$-module structure from $V(I_n(s))$.
 (1) follows. (2) follows similarly. \\
 \\
 Now consider $D(X_0, -s)$ as a dense subspace of the Hilbert representation $I_n(-s)$, which is the dual of Hilbert representation $I_n(s)$. Let $(\, , \,)$ be the pairing. 
 \commentout
 {For any $\phi \in V(I_n(s))$, $\psi
 \in D(X_0, -s)$, $u \in V(\pi)$, $v \in V(\pi)^*$, we have
 \begin{equation}
 \begin{split}
  & (\phi \otimes_{SO(p,q), D(X_0, -s) \otimes V(\pi)^*} u)(\psi \otimes v) \\
  = & \int_{SO(p,q)}(I_n(s)(g) \phi, \psi)(\pi(g)u, v) d g \\
  = & \int_{SO(p,q)}(\phi, I_n(-s)(g^{-1}) \psi)(u, \pi^*(g^{-1}) v) d g \\
  = &(\psi \otimes_{SO(p,q), V(I_n(s)) \otimes V(\pi)} v)(\phi \otimes u)
  \end{split}
  \end{equation}
  By essentially the same argument as the one preceding Definition \ref{hermitian},
  }
  By Lemma \ref{pairing},  we obtain a nondegenerate pairing 
   $$(V(I_n(s)) \otimes_{SO(p,q), D(X_0, -s) \otimes V(\pi)*} V(\pi), D(X_0, -s) \otimes_{SO(p,q), V(I_n(s)) \otimes V(\pi)} V(\pi^*)) \rightarrow \mathbb C.$$
    This pairing is $(\f o(q+d, p+d), S(O(q+d) O(p+d))$-invariant, due to the structure of $I_n(s)$ and $I_n(-s)$. Combining with (1) and (2), we obtain a nondegenerate invariant pairing between $V(I_n(s)) \otimes_{SO(p,q)} V(\pi)$ and $ D(X_0, -s) \otimes_{SO(p,q), C_c^{\infty}(X_0, s) \otimes V(\pi)} V(\pi^*)$. $\Box$. \\
  \\
  Obviously, the same statements in Theorem \ref{equivalence} hold for $-s$.  
\subsection{Subrepresentation Theorem}
By Theorem \ref{equivalence} (3) and Corollary \ref{itpkfinite}, we obtain a map
$$\mc I_1: V(I_n(s)) \otimes_{SO(p,q)} V(\pi) \rightarrow \Hom( V(\Ind_{SO(q,p) GL(d)N}^{SO(q+d, p+d)} \pi^* \otimes |\det|^{-s}), \mathbb C)$$
that preserves the $(\f o(q+d, p+d), S(O(q+d) O(p+d)))$-actions. 
Since all vectors in $V(I_n(s))$ are $S(O(q+d) O(p+d))$-finite, all vectors in $V(I_n(s)) \otimes_{SO(p,q)} V(\pi)$ are $S(O(q+d)O(p+d))$ finite.
We obtain an $(\f o(q+d, p+d), S(O(q+d) O(p+d))$ isomorphism
\begin{equation}
\begin{split}
\mc I_1: V(I_n(s)) \otimes_{SO(p,q)} V(\pi) \rightarrow & \Hom(V(\Ind_{SO(q,p) GL(d) N}^{SO(q+d, p+d)} \pi^* \otimes |\det|^{-s}), \mathbb C)_{SO(q+d)SO(p+d)} \\
\cong & V(\Ind_{SO(q,p) GL(d) N}^{SO(q+d, p+d)} \pi \otimes |\det|^s).
\end{split}
\end{equation}
Notice that the pairing between $V(I_n(s)) \otimes_{SO(p,q)} V(\pi)$ and $V(\Ind_{SO(q,p) GL(d) N}^{SO(q+d, p+d)} \pi^* \otimes |\det|^{-s})$ is nondegenerate and $S (O(q+d) O(p+d))$-invariant. Each $S (O(q+d) O(p+d))$ representation $\sigma$ has only finite multiplicity in $V(\Ind_{SO(q,p) GL(d) N}^{SO(q+d, p+d)} \pi^* \otimes |\det|^{-s})$. So the multiplicity of $\sigma^*$ in $V(I_n(s)) \otimes_{SO(p,q)} V(\pi)$ must be the same. We have shown that $\mc I_1$ is onto. We obtain
\begin{Theorem}\label{5} Let $s=\frac{n-1}{2}-m$, $n=p+q+d-1$ and $\eta=\eta(n-m-1, m)$. Let $\pi$ be an irreducible unitary representation of $SO(p,q)$ such that all leading exponent $\mu$ satisfies
$$\eta+(\mu, {\mathbf 0}_{q+d})-(d+1, d+2, \ldots, d+p, {\mathbf 0}_{q+d}) \prec 0.$$
Then $V(I_n(s)) \otimes_{SO(p,q)} V(\pi) \cong V(\Ind_{SO(q,p) GL(d) N}^{SO(q+d, p+d)} \pi \otimes |\det|^s)$.
\end{Theorem}  

\begin{cor}\label{sub}  Under the same hypothesis as Theorem \ref{5}, 
$V_{m}(-s) \otimes_{SO(p,q)} V(\pi)$ can be identified as a subrepresentation of $\Ind_{SO(q,p) GL(d) N}^{SO(q+d, p+d)} \pi \otimes |\det|^{-s}$.
\end{cor}
Proof: Recall from Theorem \ref{basicins} that $V_m(-s)$ is a submodule of $V(I_n(-s))$. Hence
$$V_{m}(-s) \otimes_{SO(p,q), V(I_n(-s))^c \otimes V(\pi)^c} V(\pi) \hookrightarrow V(I_n(-s)) \otimes_{SO(p,q)} V(\pi).$$
By Theorem \ref{5}, $V_{m}(-s) \otimes_{SO(p,q), V(I_n(-s))^c \otimes V(\pi)^c} V(\pi)$ can be identified with a Harish-Chandra submodule of  $V(\Ind_{SO(q,p) GL(d) N}^{SO(q+d, p+d)} \pi \otimes |\det|^{-s})$. \\
\\
 Let $W(-s)$ be the subspace of $V(I_n(-s))$ consisting of the $K$-types outside of
$V_{m}(-s)$.  $W(-s)$ is not $\f g$-stable and hence not a Harish-Chandra module of $SO(n,n)$. Nevertheless,
matrix coefficients $(I_n(-s)(g)V_m(-s), W(-s))$ is zero. In addition, $V(I_n(-s))=V_m(-s) \oplus W(-s)$. Hence 
$$V_{m}(-s) \otimes_{SO(p,q), V(I_n(-s))^c \otimes V(\pi)^c} V(\pi) \cong V_{m}(-s) \otimes_{SO(p,q), V_m(-s)^c \otimes V(\pi)^c} V(\pi) \cong V_{m}(-s) \otimes_{SO(p,q)} V(\pi).$$
So $V_{m}(-s) \otimes_{SO(p,q)} V(\pi)$ can be identified as a subrepresentation of $\Ind_{SO(q,p) GL(d) N}^{SO(q+d, p+d)} \pi \otimes |\det|^{-s}$.
$\Box$\\
\\
By Theorem \ref{emvm}, we have
\begin{Theorem}\label{quan_ind}
 Under the same hypothesis as Theorem \ref{5}, 
$\mc Q(2m)(V(\pi))=V(\mc E_{m}(n)) \otimes_{SO(p,q)} V(\pi)$ can be identified as a subrepresentation of $\Ind_{SO(q,p) GL(d) N}^{SO(q+d, p+d)} \pi \otimes |\det|^{m-\frac{n-1}{2}}$.
\end{Theorem}

\section{Howe's Correspondence: Nonvanishing Theorem and Minimal $K$-types}
In this section, we shall relate $\mc Q(2m)$ to the composition $\theta(2m; q+d, p+d)\theta(q,p;2m)$.
The nonvanishing of $\mc Q(2m)(V(\pi))$ then follows from Kudla's preservation principle (\cite{ku}). We will then apply Howe's theory to compute minimal $K$-types  in $\mc Q(2m)(V(\pi))$ (\cite{a} \cite{ab} \cite{mo} \cite{paul}). In the cases when a Vogan's minimal $K$-type (\cite{vogan79}) is of minimal degree in the sense of Howe (\cite{howe}), we obtain the Langlands-Vogan parameter for $\mc Q(2m)(V(\pi))$. This technique does not allow us to treat all $\mc Q(2m)(V(\pi))$ because minimal $K$-types in the sense of Vogan  may not be  of minimal degree.  \\
\\
{\bf Suppose that $2m+1 \leq n$. In this section, $\mc E_m(n)$ will be a unitary representation of $O(n,n)$, namely $\theta(2m; n, n)(\rm triv)$}.
\subsection{Associativity of Invariant Tensor Products}
Let us recall the following theorem.
\begin{Theorem}[\cite{he00}]\label{thetaanalytic}
Let $(G_1, G_2)$ be a dual reductive pair in $Sp$. Let $\omega$ be the oscillator representation of $\widetilde{Sp}$. Let $(\tilde G_1, \tilde G_2)$ be the induced covering from the metaplectic covering.
Let $\pi$ be an irreducible unitary representation of $\tilde G_1$. Then
$V(\omega) \otimes_{\tilde{G}_1} V(\pi)$, whenever defined, is equivalent to 
$ V(\theta(\pi^*))$.
\end{Theorem}
In the case of $(O(p,q), Sp_{2m}(\mb R))$, $\theta$ can always be regarded as a correspondence between representations of $O(p,q)$ and  $\widetilde{Sp}_{2m}(\mb R)$. This will be our view point from now on. 
\begin{Theorem} Let $n=p+q+d$ and $2m \leq n+1$. Let $\pi$ be an irreducible unitary representation of $O(p,q)$. Suppose that $V(\mc E_{m}(n)) \otimes_{O(p,q)} V(\pi)$ is well-defined. Then $V(\omega(q+d, p+d; 2m)) \otimes_{\widetilde{Sp}_{2m}(\mb R)} (V(\omega(p,q; 2m)) \otimes_{O(p,q)} V(\pi))$ is well-defined and as Harish-Chandra modules of $O(q+d, p+d)$
$$V(\mc E_{m}(n)) \otimes_{O(p,q)} V(\pi) \cong V(\omega(q+d, p+d; 2m)) \otimes_{\widetilde{Sp}_{2m}(\mb R)} (V(\omega(p,q; 2m)) \otimes_{O(p,q)} V(\pi)).$$
\end{Theorem}
Recall from Theorem \ref{leehe}, that 
$$V(\mc E_{m}(n)) \cong V(\omega(n,n; 2m)) \otimes_{\tilde Sp_{2m}(\mb R)} {\rm triv} \cong [V(\omega(q+d, p+d; 2m)) \otimes V(\omega(p,q; 2m)] \otimes_{{\widetilde Sp}_{2m}(\mb R)} {\rm triv}.$$
So when we restrict our representation onto $O(p,q) O(q+d, p+d)$, we have
$$V(\mc E_{m}(n))|_{O(p,q) O(q+d, p+d)} \cong V(\omega(p,q; 2m)) \otimes_{{\widetilde Sp}_{2m}(\mb R)} V(\omega(q+d, p+d; 2m)).$$
Then it suffices to prove the associativity of invariant tensor products
\begin{equation}
\begin{split}
 & V(\omega(q+d, p+d; 2m)) \otimes_{\widetilde{Sp}_{2m}(\mb R)} [V(\omega(p,q; 2m)) \otimes_{O(p,q)} V(\pi)]\\
  \cong & [V(\omega(q+d, p+d; 2m)) \otimes_{\widetilde{Sp}_{2m}(\mb R)} V(\omega(p,q; 2m)) ] \otimes_{O(p,q)} V(\pi).
  \end{split}
  \end{equation}
Theorem \ref{xios} and Lemma \ref{kfinitees} allow us to apply Fubini's theorem.  Then we can  interchange  two integrals defined by  the averaging operator $\mc A_{\widetilde{Sp}_{2m}(\mb R)}$ and $\mc A_{O(p,q)}$.\\
  \\
Proof: Suppose that $V(\mc E_{m}(n)) \otimes_{O(p,q)} V(\pi)$ is well-defined.
\begin{enumerate}
\item 
Let $\exp -\frac{1}{2} \|x \|^2$ and $\exp -\frac{1}{2} \| y \|^2$ be the lowest weight vectors in the Schr\"odinger model of $\omega(p,q;2m)$ and $\omega(q+d,p+d; 2m)$. Then 
$(\exp -\frac{1}{2} \|x \|^2) \otimes_{\widetilde{Sp}_{2m}(\mb R)} (\exp -\frac{1}{2} \| y \|^2 )$ yields a nonzero spherical vector in $\mc E_m(n)$.  By Theorem \ref{xios}, $\forall$ $ g_1 \in O(p,q)$, $g_2 \in O(q+d, p+d)$,
\begin{equation}
\begin{split}
\Xi_{\eta}(g_1, g_2)=C \int_{h \in \widetilde{Sp}_{2m}(\mb R)} & (\omega(p,q;2m)(g_1, h)\exp -\frac{1}{2} \|x \|^2, \exp -\frac{1}{2} \|x \|^2 ) \\
& (\omega(q+d,p+d; 2m)(g_2, h) \exp -\frac{1}{2} \| y \|^2, \exp -\frac{1}{2} \| y \|^2) d h.
\end{split}
\end{equation}
and the integrand is always positive. 
Since  the vector $(\exp -\frac{1}{2} \|x \|^2 \otimes_{\widetilde{Sp}_{2m}(\mb R)} \exp -\frac{1}{2} \| y \|^2 )\otimes_{O(p,q)} v$ is well-defined for all $v \in V(\pi)$, we have
$\Xi_{\eta}(g_1) |(\pi(g_1) u, v)| \in L^1(O(p,q))$. By Fubini's theorem,
$$ (\omega(p,q;2m)(g_1, h)\exp -\frac{1}{2} \|x \|^2, \exp -\frac{1}{2} \|x \|^2 )$$
$$(\omega(q+d,p+d; 2m)(h) \exp -\frac{1}{2} \| y \|^2, \exp -\frac{1}{2} \| y \|^2) |(\pi(g_1) u, v)|$$
is integrable on $(g_1, h) \in O(p,q) \times \widetilde{Sp}_{2m}(\mb R)$. 
\item In particular, for almost all $h \in \widetilde{Sp}_{2m}(\mb R)$, we have
\begin{equation}\label{omegapq2mlowest}
|(\omega(p,q;2m)(g_1, h)\exp -\frac{1}{2} \|x \|^2, \exp -\frac{1}{2} \|x \|^2 )(\pi(g_1) u, v)| \in L^1(O(p,q)).
\end{equation}
From Lemma \ref{basices}, we have the exact form
\begin{equation}
\begin{split}
  & |(\omega(p,q;2m)(g_1, h)\exp -\frac{1}{2} \|x \|^2, \exp -\frac{1}{2} \|x \|^2 )| \\
  = & [\prod_{i=1}^p \prod_{j=1}^m (a_i b_j+ a_i^{-1} b_j^{-1})^{-\frac{1}{2}} (a_i b_j^{-1}+a_i^{-1} b_j)^{-\frac{1}{2}}] [\prod_{j=1}^m (b_j+b_j^{-1})^{-\frac{q-p}{2}}] \\
= & [\prod_{i=1}^p \prod_{j=1}^m (a_i^2+a_i^{-2}+b_j^2+b_j^{-2})^{-\frac{1}{2}}] [\prod_{j=1}^m (b_j+b_j^{-1})^{-\frac{q-p}{2}}],
\end{split}
\end{equation}
where $g_1=k_1 a k_2$ and $h=u_1 b u_2$ are the $KA^+K$ decompositions for $O(p,q)$ and $\widetilde{Sp}_{2m}(\mb R)$. See \cite{basic} for more details. 
The integrability of (\ref{omegapq2mlowest}) for one $h \in \widetilde{Sp}_{2m}(\mb R)$ implies the integrability of $(\ref{omegapq2mlowest})$ for the identity element $e$. This is essentially due to the fact that
for fixed $b$,
$$\frac{1}{a_i^2+a_i^{-2}+b_j^2+b_j^{-2}} \geq \frac{2}{(b_j^2+b_j^{-2})(a_i^2+a_i^{-2}+2)}.$$
It follows that $V(\omega(p,q; 2m)) \otimes_{O(p,q)} V(\pi)$ is well-defined. By Theorem \ref{thetaanalytic}, 
$$V(\omega(p,q; 2m)) \otimes_{O(p,q)} V(\pi)=V(\theta(p,q;2m)(\pi^*)).$$ 
\item Let ${}^{\widetilde{Sp}_{2m}(\mb R)} V(\omega(p,q;2m))$ be the finite linear span of $\widetilde{Sp}_{2m}(\mathbb R)$ translations of $V(\omega(p,q;2m)$. Indeed $ {}^{\widetilde{Sp}_{2m}(\mb R)} V(\omega(p,q;2m)) \otimes_{O(p,q)} V(\pi)$ is well-defined. This can be seen as follows.
Notice
$$\frac{1}{a_i^2+a_i^{-2}+b_j^2+b_j^{-2}} \leq \frac{1}{(a_i+a_i^{-1})(b_j+b_j^{-1})},$$
the integrability of $(\ref{omegapq2mlowest})$ at $h=e$ then implies the integrability of $(\ref{omegapq2mlowest})$ for every $h \in \widetilde{Sp}_{2m}(\mb R)$. Here $b$ is the $A^+$ component of the $KA^+K$ decomposition of $h$.
By Lemma \ref{kfinitees}, for any $\phi_1, \phi_2 \in \mc P(x), h \in \widetilde{Sp}_{2m}(\mb R)$,
$$|(\omega(p,q;2m)(g_1, h) (\phi_1(x) \exp -\frac{1}{2} \|x \|^2), \phi_2 \exp -\frac{1}{2} \|x \|^2 )  |(\pi(g_1) u, v)| \in L^1(O(p,q)).$$
\item Now $V(\omega(p,q; 2m)) \otimes_{O(p,q)} V(\pi)$ is a Harish-Chandra module for a continuous representation $\theta(p,q;2m)(\pi^*)$.  Since $\omega(p,q;2m)$ and $\pi$ are both unitary, by Prop. \!\ref{invher},
the canonical Hermitian form is $( \f{sp}_{2m}(\mb R), \tilde U(m))$-invariant.
By Lemma \ref{kfinitees} and Part (1) of the proof, we have for any $\psi_1, \psi_2 \in \mc P(y)$
\begin{equation}
\begin{split}
 & (\omega(p,q;2m)(g_1, h) (\phi_1(x) \exp -\frac{1}{2} \|x \|^2), \phi_2 \exp -\frac{1}{2} \|x \|^2 ) \\
 & (\omega(q+d,p+d, 2m)(h) (\psi_1(y) \exp -\frac{1}{2} \| y \|^2), \psi_2 (y) \exp -\frac{1}{2} \| y \|^2) |(\pi(g_1) u, v)|
 \end{split}
 \end{equation}
is integrable on $ O(p,q) \times \widetilde{Sp}_{2m}(\mb R)$. 
 So $V(\omega(q+d, p+d; 2m)) \otimes_{\widetilde{Sp}_{2m}(\mb R)} [V(\omega(p,q; 2m)) \otimes_{O(p,q)} V(\pi)]$ is well-defined.
 It can be constructed using double integral on $(g_1, h) \in O(p,q) \times \widetilde{Sp}_{2m}(\mb R)$.
\\
\\
 By Fubini's theorem, we have 
 $$V(\mc E_{m}(n)) \otimes_{O(p,q)} V(\pi) \cong V(\omega(q+d, p+d; 2m)) \otimes_{\widetilde{Sp}_{2m}(\mb R)} (V(\omega(p,q; 2m)) \otimes_{O(p,q)} V(\pi)).$$
 In addition, both spaces inherit the same $(\f o(q+d, p+d), O(q+d)O(p+d))$-module structure from $V(\omega(q+d+p, p+d+q; 2m))$. Our theorem follows immediately.
 \end{enumerate}
 $\Box$

\subsection{Nonvanishing of Quantum Induction}
 Identify $O(p,q)$ with $O(q,p)$ as in Lemma \ref{basepoint}. This amounts to essentially rearranging the coordinates. A representation $\pi$ of $SO(p,q)$ will also be  regarded as a representation of $SO(q, p)$.
We would like to relate $\mc Q(2m)(V(\pi))$ to $\theta(2m; q+d, p+d) \theta(q,p; 2m) (\Ind_{SO(q,p)}^{O(q,p)} \pi)$.
\begin{lem}\label{pqqp} As unitary representations of $(O(p,q), \widetilde{Sp}_{2m}(\mb R))$,  $\omega(p,q;2m) \cong \omega(q, p; 2m)^* \cong \omega(q, p; 2m)^c$. In addition, $\theta(p,q; 2m)(\pi^c) \cong [\theta(q,p; 2m)(\pi)]^c$.
\end{lem}
Proof: Let $\omega_m$ be the oscillator representation of $\widetilde{Sp}_{2m}(\mb R)$. Then $\omega(p,q;2m)$ can be modelled by $[\otimes^p \omega_m] \otimes [\otimes^q \omega_m^c]$. So
$\omega(p,q;2m)^c$ can be modelled by $[\otimes^p \omega_m^c] \otimes [\otimes^q \omega_m]$.
By reordering of coordinates, we have  $\omega(p,q; 2m)^c \cong \omega(q, p; 2m)$. \\
\\
If $\pi \otimes \theta(p,q;2m)(\pi)$ appears as a quotient of $\omega(p,q; 2m)$, then
$\pi^c \otimes [\theta(p,q;2m)(\pi)]^c$ appears as  a quotient of $\omega(p,q;2m)^c \cong \omega(q,p;2m)$.
So $\theta(q,p;2m)(\pi^c) \cong \theta(p,q; 2m)(\pi)^c$. Equivalently,
$$\theta(p,q;2m)(\pi^c) \cong [\theta(q,p; 2m)(\pi)]^c.$$
$\Box$ 
 \begin{Theorem}
 If $2m+1 \geq p+q$, then either $\theta( 2m; q+d,p+d)(\theta(q, p; 2m)(\pi)) \neq 0$ or $\theta( 2m; q+d,p+d)(\theta(q, p; 2m)(\pi \otimes \det)) \neq 0$.
 \end{Theorem}
  Proof:
By a theorem of Moeglin (\cite{mo}) and a theorem of Adams-Barbasch (\cite{ab}), we have either $\theta(q, p; 2m)(\pi) \neq 0$ or $
\theta(q, p; 2m)(\pi \otimes \det) \neq 0$. See also \cite{non}. By Kudla's preservation principle (\cite{ku}), $\theta( 2m; q+d,p+d)(\theta(q, p; 2m)(\pi)) \neq 0$ if $\theta(q, p; 2m)(\pi)) \neq 0$. Our theorem then follows. $\Box$ \\
\\
Now suppose that $\overline{\pi}$ is an irreducible unitary representation of $O(p,q)$, then ${\overline{\pi}}^* \cong {\overline{\pi}}^c$. By a Theorem of Przebinda \cite{prz}, $\theta({\overline{\pi}})$ will have an invariant Hermitian structure, i.e., 
$[\theta(\cdot)({\overline{\pi}})]^* \cong [\theta(\cdot)({\overline{\pi}})]^c$. We obtain
\begin{equation}\label{quanandtheta}
\begin{split}
& V(\mc E_m(n)) \otimes_{O(p,q)} V({\overline{\pi}}) \\
 \cong & V(\omega(q+d, p+d; 2m)) \otimes_{\widetilde{Sp}_{2m}(\mb R)} (V(\omega(p,q; 2m)) \otimes_{O(p,q)} V({\overline{\pi}}))
 \\
\cong & V(\omega(q+d, p+d; 2m)) \otimes_{\widetilde{Sp}_{2m}(\mb R)} V(\theta(p,q;2m)({\overline{\pi}}^*))  \commentout{(\mathrm{Theorem \,\, \ref{thetaanalytic}})} \\
  \cong & V(\theta(2m; q+d,p+d)( \theta(p,q; 2m)({\overline{\pi}}^*)^*)) \commentout{( \mathrm{Theorem \,\, \ref{thetaanalytic}} )}\\
 \cong & V(\theta(2m; q+d,p+d)( \theta(p,q; 2m)({\overline{\pi}}^c)^c)) \commentout{( \mathrm{Theorem \,\,  5.10  \,\, \cite{prz}}) }\\
\cong & V(\theta(2m; q+d,p+d)(\theta(q, p; 2m)(\overline{\pi})))  \commentout{( \mathrm{Lemma} \,\, \ref{pqqp})}.
\end{split}
 \end{equation}

\begin{Theorem}\label{quan_theta} Suppose that $2m+1 \geq p+q$. Let $\pi$ be an irreducible unitary representation of $SO(p,q)$. Then $\mc Q(2m)(V(\pi)) \neq 0$. As Harish-Chandra modules of $SO(q+d, p+d)$, $\mc Q(2m)(V(\pi))$ is equivalent to $\theta( 2m; q+d,p+d)(\theta(q, p; 2m)(\overline{\pi}))$ if $\overline{\pi}$ is unique, is equivalent to 
$$\theta( 2m; q+d,p+d)(\theta(q, p; 2m)(\overline{\pi})) \oplus \theta( 2m; q+d,p+d)(\theta(q, p; 2m)(\overline{\pi} \otimes \det))$$
if $\overline{\pi}$ is not unique.
\end{Theorem}
Proof:  First of all, we have
\begin{equation}
\begin{split}
V(\mc E_m(n)) \otimes_{SO(p,q)} V(\pi) \cong  & V(\mc E_m(n)) \otimes_{O(p,q)}(V(\Ind_{SO(p,q)}^{O(p,q)} \pi)). \\
\end{split}
\end{equation}
If 
$\Ind_{SO(p,q)}^{O(p,q)} \pi$ is irreducible, then $\overline{\pi} \cong \Ind_{SO(p,q)}^{O(p,q)} \pi$ and 
$ (\Ind_{SO(p,q)}^{O(p,q)} \pi) \otimes \det \cong \Ind_{SO(p,q)}^{O(p,q)} \pi.$ By the previous theorem and Eq. \ref{quanandtheta}, $$\mc Q(2m)(V(\pi)) \cong \theta( 2m; q+d,p+d)\theta(q, p; 2m)(\overline{\pi}) \neq 0.$$
If 
$\Ind_{SO(p,q)}^{O(p,q)} \pi$ is not irreducible, then it must contain two subrepresentations $\overline{\pi}$ and $\overline{\pi} \otimes \det$, where $\overline{\pi}$ is a unitary representation of $O(p,q)$. By the previous theorem, $$\mc Q(2m)(V(\pi)) \cong \theta( 2m; q+d,p+d)(\theta(q, p; 2m)(\overline{\pi})) \oplus \theta( 2m; q+d,p+d)(\theta(q, p; 2m)(\overline{\pi} \otimes \det)) \neq 0.$$
$\Box$
\subsection{Langlands-Vogan Parameter of Quantum Induced Module}
Howe's correspondence uniquely determines a one to one correspondence 
between the $K$-types of the lowest degrees in $\pi$ and $\theta(\pi)$ (\cite{howe}).
{\bf We denote this $K$-types correspondence by $\theta_0$. For any $K$-type $\tau$, let $d(\tau, \omega)$ to be the lowest degree of $\tau$-isotypic subspace in the Fock-Segal-Bargman model of $\omega$.} $\theta_0$ is often called the correspondence of joint harmonics.  $\theta_0$ is known explicitly (\cite{kv} \cite{mo} \cite{a}). In the instances that a minimal
$K$-type in the sense of Vogan is a $K$-type of lowest degree, Howe's correspondence can be determined in terms of the Langlands-Vogan parameter. This has been done mostly in the equal rank case (\cite{mo},\cite{ab}, \cite{paul}).
With the appearance of $(\xi, -)$ for $O(p)$ or $O(q)$, a minimal $K$-type may not be a $K$-type of the lowest degree  and Howe's correspondence is difficult to compute explicitly. \\
\\
Two necessary conditions for defining nonvanishing quantum induction in Theorems \ref{quan_ind} and \ref{quan_theta} are $p+q \leq 2m+1$ and $s=  \frac{p+q+d-1}{2}-m \geq 0$. {\bf We suppose that $p+q \leq 2m+1 \leq p+q+d$. Suppose that $((\xi_1, +),(\xi_2, +))$ is a minimal $O(p) O(q)$ type of an irreducible constituent $\overline{\pi}$ of $\Ind_{SO(p,q)}^{O(p,q)} \pi$}.  The minimal degree of $((\xi_2, +), (\xi_1, +))$ in $\omega(q,p;2m)$  is  the sum of all entries of $\xi_1$ and $\xi_2$ (\cite{mo} \cite{ab} \cite{a}). The key observation is that 
$d(((\xi_2, +), (\xi_1, +)), \omega(q,p;2m))$ is independent of $m$ as long as $p+q \leq 2m+1$. By the results of Moeglin, Adam-Barbasch and Paul, for the equal size cases,
$$d(((\xi_2, +), (\xi_1, +)), \omega(q,p;2m))=\min \{ d(\tau, \omega(q,p;2m)) \mid \tau \subseteq V(\overline{\pi}), \tau \in \widehat{O(q)O(p)} \},$$
i.e.,  the minimal $K$-type $(\xi_2, +) \otimes (\xi_1, +)$ is a $K$-type of minimal degree. It follows that for any $2m \geq p+q-1$, we have
$$d(((\xi_2, +), (\xi_1, +)), \omega(q,p;2m))=\min \{ d(\tau, \omega(q,p;2m)) \mid \tau \subseteq V(\overline{\pi}), \tau \in \widehat{O(q)O(p)} \}.$$
Now $\theta(q,p; 2m)(\ol{\pi})$ is an irreducible representation of $\widetilde{Sp}_{2m}(\mb R)$. By Howe's theory (\cite{howe}),  the corresponding $K$-type 
$$\theta_0(q,p;2m)((\xi_2, +) \otimes (\xi_1, +))=(\xi_2, {\mathbf 0}, -\xi_1)+{\bf \frac{q-p}{2}}.$$ This $K$-type in $\theta(q,p;2m)(\ol{\pi})$
 has a similar property:
 $$d(((\xi_2, {\mathbf 0}, -\xi_1)+{\bf \frac{q-p}{2}}), \omega(q,p;2m))=\min \{ d(\tau^{\prime}, \omega(q,p;2m)) \mid \tau^{\prime} \subseteq V(\theta(q,p;2m)(\ol{\pi})) \}.$$
 Observe that $ d(\tau^{\prime}, \omega(q,p;2m)) = d (\tau^{\prime}, \omega(q+d, p+d; 2m))$. So the statement above holds if we replace $\omega(q,p;2m)$ by $\omega(q+d,p+d;2m)$. 
 It follows that the $K$-type $(\xi_2, {\mathbf 0}, -\xi_1)+{\bf \frac{q-p}{2}}$ of $\theta(q,p;2m)(\ol{\pi})$ is also of the minimal degree in $\omega(2m; q+d, p+d)$. Hence 
 $$\theta_0(2m; q+d, p+d)({(\xi_2, {\mathbf 0}, -\xi_1)}+{\bf \frac{q-p}{2}}) \subseteq V(\theta(2m; q+d, p+d) \theta(q, p; 2m)(\ol{\pi})).$$
The left hand side is equal to ${((\xi_2 \oplus {\mathbf 0}, +),(\xi_1 \oplus {\mathbf 0}, +))}$. Hence the $K$-type 
$((\xi_2 \oplus {\mathbf 0}, +),(\xi_1 \oplus {\mathbf 0}, +))$ occurs in $\theta(2m; q+d, p+d) \theta(q, p; 2m)(\ol{\pi})$.\\
\\
By Theorem \ref{quan_theta}, $\mc Q_{2m}(V(\pi))$ will contain the $K$-type $(\xi_2 \oplus {\mathbf 0}, \xi_1 \oplus {\mathbf 0}, +)$. 
By Theorem \ref{quan_ind}, $\mc Q_{2m}(V(\pi))$ is a quotient of $\Ind_{SO(q,p) GL(d)N }^{SO(q+d, p+d)} \pi \otimes |\det|^s$.
Notice that $((\xi_2 \oplus {\mathbf 0}), (\xi_1 \oplus {\mathbf 0}), +)$ is the  minimal $K$-type of $\Ind_{SO(q,p) GL(d)N }^{SO(q+d, p+d)} \pi \otimes |\det|^s$ (Ch 8. \cite{gw}).
\begin{Theorem}\label{quan_vogan}
Suppose that $p+q \leq 2m+1 \leq p+q+d$. Suppose that $(\xi_1, \xi_2, +)$ is a minimal $S(O(p) O(q))$ type of ${\pi}$. Then $\mc Q_{2m}(V(\pi))$ contains the irreducible Vogan subquotient of $\Ind_{SO(q,p) GL(d)N }^{SO(q+d, p+d)} \pi \otimes |\det|^s$.
Here $s=\frac{n-1}{2}-m=\frac{p+q+d-2m-1}{2}$.
\end{Theorem}

 \section{Unitary Langlands-Vogan Parameters  }
 The main purpose of this section is to determine the unitarity of certain Langlands-Vogan parameter. By a theorem of Harish-Chandra, an admissible representation $\pi$ is unitarizable if $V(\pi)$ has a $(\f g, K)$-invariant pre-Hilbert structure. By Prop. \ref{invher},
 $\mc Q(2m)(V(\pi))$  has an invariant Hermitian form.  Theorem \ref{posinv} allows us to determine when this Hermitian form will be positive definite. 
We have the following
\begin{Theorem}[Theorem A] 
Suppose that $p+q \leq 2m+1 \leq p+q+d$ and $p \leq q$. Let $\pi$ be an irreducible unitary representation of $SO(p,q)$ such that
its every $K$ finite matrix coefficient $f(g)$ satisfies the condition that
$$ |f(g)| \leq C_f \exp (\overbrace{m+2-p-q- \epsilon, m+3-p-q, \ldots m+1-q}^{p})(| H^+(g)|).$$
for some $\epsilon >0$. Then
\begin{enumerate}
\item $\mc Q(p,q; 2m; q+d, p+d)(V(\pi))$ is well-defined;
\item $\mc Q(p,q; 2m; q+d, p+d)(V(\pi))$ is nonvanishing;
\item $\mc Q(p,q; 2m; q+d, p+d)(V(\pi))$ is a subrepresentation of
$\Ind_{SO(q,p) GL(d) N}^{SO(q+d, p+d)} \pi \otimes \mathbb |\det|^{m-\frac{p+q+d-1}{2}}$ upon identifying $SO(p,q)$ with $SO(q,p)$;
\item $\mc Q(p,q; 2m; q+d, p+d)(V(\pi))$ is unitarizable.
\end{enumerate}
\end{Theorem}
Proof: (1) follows from Theorem \ref{posinv}. (2) follows from Theorem \ref{quan_theta}. (3) follows from Cor \ref{sub}. We only need to prove $(4)$. Notice that 
$$V(\mc E_m(n)) \otimes_{SO(p,q)} V(\pi)=V(\mc E_m(n)) \otimes_{O(p,q)} V(\Ind_{SO(p,q)}^{O(p,q)} \pi).$$
is either irreducible or decompose into two irreducible admissible representations of $O(q+d, p+d)$. Each irreducible one will be of the form $V(\mc E_m(n)) \otimes_{O(p,q)} u$ for some $u \in V(\Ind_{SO(p,q)}^{O(p,q)} \pi)$. Applying Prop. \ref{posinv} to $O(p,q)$, the invariant Hermitian form on $V(\mc E_m(n)) \otimes_{O(p,q)} u$ is positive definite.  By Harish-Chandra's theorem, $\mc Q(p,q; 2m; q+d, p+d)(V(\pi))$ is unitarizable. $\Box$

\begin{Theorem}[Theorem B]
Suppose that $p+q \leq 2m+1 \leq p+q+d$ and $p \leq q$. Let $\pi$ be an irreducible unitary representation of $SO(q,p)$ such that
its every $K$ finite matrix coefficient $f(g)$ satisfies the condition that
$$ |f(g)| \leq C_f \exp (\overbrace{m+2-p-q- \epsilon, m+3-p-q, \ldots m+1-q}^{p})(| H^+(g)|).$$
for some $\epsilon >0$. 
Suppose that the minimal $S(O(q) O(p))$-types of $ \pi$ contain
a $(\xi, \eta, +)$. Then the Vogan subquotient of
$$\Ind_{SO(q,p) GL(1)^{d} N}^{SO(q+d, p+d)} \pi \otimes |\det|^{m-\frac{p+q}{2}} |\det|^{m-1-\frac{p+q}{2}} \ldots |\det|^{m-\frac{p+q}{2}-d+1}$$
is unitary.
\end{Theorem}
Proof: Since
$${\bf {m-\frac{p+q+d-1}{2}}}-\rho(SL_d(\mb R))=(m-\frac{p+q}{2}-d+1, m-\frac{p+q}{2}-d+2, \ldots, m-\frac{p+q}{2}),$$
the Vogan subquotient of $\Ind_{SO(q,p) GL(d)N }^{SO(q+d, p+d)} \pi \otimes |\det|^{\frac{2m-p-q-d+1}{2}}$ is exactly the Vogan subquotient of 
$$\Ind_{SO(q,p) GL(1)^{d} N}^{SO(q+d, p+d)} \pi \otimes |\det|^{m-\frac{p+q}{2}} |\det|^{m-1-\frac{p+q}{2}} \ldots |\det|^{m-\frac{p+q}{2}-d+1}.$$
Our theorem follows from Theorem \ref{quan_vogan} and Theorem [A] (4). $\Box$

\subsection{Arthur's Packet}
Let $G$ be an inner form of an algebraic reductive group. Let ${}^L G$ be the Langlands dual group. 
Let $W_{\mathbb R}= \mathbb Z_2  \ltimes \mathbb C^{\times}$ be the Weil group.
A Langlands parameter $\phi: W_{\mathbb R} \rightarrow {}^L G$ can be decomposed into a product 
of a compact (also known as tempered) part $\phi_0$ and a noncompact part $\phi_{+}$. $\phi_0$ and $\phi_{+}$ commute with each other. Arthur's parameter is a map
$$\psi: W_{\mathbb R} \times SL(2, \mb C) \rightarrow {}^L G,$$
such that $\psi|_{W_{\mathbb R}}$ is a tempered parameter.
Since $\psi(W_{\mathbb R})$ and $\psi(SL(2, \mb C))$ commute,
$\psi(SL(2, \mathbb C))$ must be in the centralizer of $\psi(W_{\mathbb R})$, in fact, the identity component of the centralizer of  $\psi(W_{\mathbb R})$. Let $C(\psi(W_{\mb R}))_0$ be the identity component of the centralizer. $C(\psi(W_{\mb R}))_0$ is a reductive group. Then Arthur defines a Langlands parameter $$\phi_{\psi}: \phi_{\psi}(w)=\psi(w, H_{|w|}), (w \in W_{\mathbb R})$$
where $H_{|w|}={\rm diag}(|w|^{\frac{1}{2}}, |w|^{-\frac{1}{2}}) \in SL(2, \mathbb C)$. Arthur conjectured that the representations in $\phi_{\psi}$ are unitary (\cite{ar}). \\
\subsection{Type $D$: $p+q$ even}
In the case $G=SO(p,q)$ with $p+q$ even, ${}^L G_0= SO(p+q, \mb C)$. We consider only those $\psi$ with
$\psi(w)$ a tempered representation of $SO(p-k, q-k)$ and $\psi(SL(2, \mathbb C)) \subseteq SO(2k, \mb C)$. We recall the following facts (\cite{cm}).
\begin{enumerate}
\item Group homomorphisms from $SL(2, \mb C)$ to $SO(2k, \mathbb C)$ are in one-to-one correspondence to Lie algebra homomorphisms from the standard triple $\{ H, X, Y \}$ to the Lie algebra $\f{so}(2n, \mb C)$.
\item The $SO(2n, \mb C)$-conjugacy classes of Lie algebra homomorphisms from the standard triple $\{ H, X, Y \}$ to the Lie algebra $\f{so}(2n, \mb C)$ are in one-to-one correspondence with nilpotent adjoint orbits of $SO(2n, \mb C)$, namely the $SO(2n, \mb C)$ conjugacy class of $X$. 
\item Nilpotent adjoint orbits of $O(2k, \mathbb C)$ are in one-to-one correspondence with orthogonal Young diagrams of size $2k$. A Young diagram is called orthogonal if each even part occurs with even multiplicity.
\end{enumerate}
Denote the nilpotent adjoint orbit of $O(2k, \mathbb C)$ corresponding to $\mathbf D$ by $\mc O_{\mathbf D}$. We say $\mathbf D$ is very even if only even parts occur in $\mathbf D$. 
Almost all nilpotent adjoint orbits $\mc O_{\mathbf D}$ of $O(2k, \mathbb C)$ are nilpotent adjoint orbits of $SO(2k, \mathbb C)$, except for  very even $\mathbf D$. For very even $\mathbf D$, $\mc O_{\mathbf D}$  splits into two nilpotent orbits of $SO(2k, \mathbb C)$. We denote them by $\mc O_{\mathbf D, \pm}$.  \\
\\
Let $\mathbf D$ be an orthogonal Young diagram, given by the partition of $2k$:
$$\sum_{i=1}^r d_i=2k, \qquad (d_1 \leq d_2 \leq d_3 \ldots \leq d_r ) .$$
In order to state Arthur's conjecture, we must know $s_{\mathbf D}$, the semisimple element in $\mathfrak{so}(2k, \mathbb C)$ corresponding to ${\rm diag}(\frac{1}{2}, -\frac{1}{2})$ in $\f{sl}(2, \mb C)$. Parametrize
 the conjugacy class of $s_{\mathbf D}$ by $v_{\mathbf D}$. $v_{\mathbf D}$ shall be identified with an element in $\mathbb R^k$ up to the action of the Weyl group of $SO(2k, \mathbb C)$. $v_{\mathbf D}$ can be constructed as a direct sum in the following way (Pg. 78-79 \cite{cm}).
 \begin{enumerate}
 \item If $d_i=d_{i+1}=d$, define $v_{(d,d)}=(\frac{d-1}{2}, \frac{d-3}{2}, \ldots, \frac{1-d}{2})$. We construct $v_{(d,d)}$ for all pairs. Since even $d$ appears with even multiplicities, we will only be left with odd $d$'s. There are even number of them.
 \item If $(d_i, d_j)$ is an odd pair, define $v_{(d_i, d_j)}=(\frac{d_i-1}{2}, \frac{d_i-1}{2}-1, \ldots 1, 0, -1, \ldots \frac{1-d_j}{2})$. Due to the Weyl group action, $v_{(d_i, d_j)}$ is equivalent to $(\frac{d_i-1}{2}, \frac{d_i-1}{2}-1, \ldots, 1, 0, 1, 2, \ldots \frac{d_j-1}{2})$.
 \item 
 If $\mathbf D$ is not very even, there must be at least two odd $d_i$'s. Then $0$ appears in $v_{\mathbf D}$. Due to the Weyl group action, we can use $|v_{\mathbf D}|$ to represent the semisimple element $s_{\mathbf D}$. If $\mathbf D$ is very even, we can apply the Weyl group action to make the entries of $v_{\mathbf D}$ all positive, except possibly one entry. Then there are two nonequivalent $v_{\mathbf D, \pm}$;
 $v_{\mathbf D, +}$ with even number of negative numbers and $v_{\mathbf D, -}$ with odd number of negative numbers.
 Obviously, $\mc O_{\mathbf D, \pm}$ are related to each other by the group action $O(2k, \mathbb C)/SO(2k, \mathbb C)$.
 \end{enumerate}
 Now we can restate Theorem B for $p+q$ even.
 \begin{cor}\label{typeb1}
 Suppose that $p+q \leq 2m+1 \leq p+q+d$, $p \leq q$ and $p+q$ even. Let $\pi$ be an irreducible unitary representation of $SO(p,q)$ such that
its every $K$ finite matrix coefficient $f(g)$ satisfies the condition that
$$ |f(g)| \leq C_f \exp (\overbrace{m+2-p-q- \epsilon, m+3-p-q, \ldots m+1-q}^{p})(| H^+(g)|).$$
for some $\epsilon >0$. 
Suppose that the minimal $S(O(p) O(q))$-types  contain
a $(\xi, \eta, +)$. Then Vogan subquotient of 
$\Ind_{SO(p,q) GL(1)^{d} N}^{SO(p+d, q+d)} \pi \otimes ({\rm triv} \otimes \mathbb C_{ v_{(2m+1-p-q, p+q+2d-2m-1)}})$
is unitary.
 \end{cor}
 \begin{re}\label{typeb2} Now let $2m+1-p-q=s$ and $p+q+2d-2m-1=t$. Then $s,t$ are odd. The growth condition becomes
 $$ |f(g)| \leq C_f \exp (\overbrace{-\frac{p+q}{2}+\frac{s+3}{2}- \epsilon, -\frac{p+q}{2}+\frac{s+5}{2}, \ldots \frac{-q+p}{2}+\frac{s+1}{2}}^{p})(| H^+(g)|).$$
 The assertion is that the Vogan subquotient of $\Ind_{SO(p,q) GL(1)^{d} N}^{SO(p+d, q+d)} \pi \otimes ({\rm triv} \otimes \mathbb C_{ v_{(s,t)}})$ is unitary.
 \end{re}
 \begin{Theorem}[Theorem C] Let $p+q$ be even. Let $\mathbf D$ be an orthogonal Young diagram of size $2k$. Let $\sigma$ be an irreducible tempered representation of $SO(p-k, q-k)$ with a minimal $S(O(p-k)O(q-k))$-type $(\xi, \eta, +)$. Then the Langlands-Vogan parameter $(SO(p-k, q-k)GL(1)^k N, \sigma \otimes {\rm triv}, v_{\mathbf D})$ is unitary, with a minimal $K$-type
 $(\xi \oplus {\mathbf 0}, \eta \oplus {\mathbf 0}, +)$.
 \end{Theorem}
 Proof: Clearly, a minimal $K$-type of $\Ind_{SO(p-k, q-k)GL(1)^k N}^{SO(p,q)} \sigma \otimes ({\rm triv} \otimes  \mathbb C_{ v_{\mathbf D}})$ is $(\xi \oplus {\mathbf 0}, \eta \oplus {\mathbf 0}, +)$. Let $\mathbf D$ be defined by the partition
 $$d_1 \leq d_2 \leq d_3 \ldots \leq d_r \qquad (\sum_{i=1}^r d_i=2k).$$ 
 We apply induction on $r$. If $r=0$, we have the tempered representation $\sigma$ which is unitarizable.\\
 \\
 If there is a pair $d_i=d_{i+1}=d$, let $\mathbf D_0$ be the Young diagram obtained by deleting $d_i, d_{i+1}$. Then by induction hypothesis, the Langlands-Vogan parameter $(SO(p-k, q-k)GL(k-d)N, \sigma \otimes {\rm triv}, v_{\mathbf D_0})$ is unitary. Denote this representation by $\pi_0$. Then $\Ind_{SO(p-d, q-d) GL(d) N}^{SO(p,q)} \pi_0 \otimes {\rm triv}$ is unitary. It is easy to see that the Langlands-Vogan parameter $(SO(p-k,q-k) GL(1)^k N, \sigma \otimes {\rm triv}, v_{\mathbf D})$ is a subrepresentation of 
 $\Ind_{SO(p-d, q-d) GL(d) N}^{SO(p,q)} \pi_0 \otimes {\rm triv}$, hence is unitary. \\
 \\
 If there are no repeats in $d_i$, then $d_1, d_2, \ldots, d_r$ must be all odd and $r$ must be even.
 If $r=2$, set $d_0=0$.
 Let $\mathbf D_0$ be the Young diagram obtained from $\mathbf D$ by deleting $d_r$ and $d_{r-1}$.  
 By induction hypothesis, the Langlands-Vogan parameter $(SO(p-k,q-k) GL(1)^{k-\frac{d_r+d_{r-1}}{2}} N, \sigma \otimes {\rm triv}, v_{\mathbf D_0})$ is unitary. Denote this unitary representation of $SO(p-\frac{d_r+d_{r-1}}{2}, q-\frac{d_r+d_{r-1}}{2})$ by $\pi_0$. Then by  Remark \ref{langlandsmatrixco}, $K$-finite matrix coefficients of $\pi_0$ are bounded by multiples of $\Xi_{v_{\mathbf D_0} \oplus {\bf 0_{p-k}}}(g)$. The leading exponents of $\Xi_{v_{\mathbf D_0} \oplus {\bf 0_{p-k}}}(g)$ is bounded by 
 \begin{equation}
 \begin{split}
  & (\frac{d_{r-2}}{2}, \frac{d_{r-2}}{2}, \ldots, \frac{d_{r-2}}{2}) -(\frac{p+q-d_r-d_{r-1}}{2}-1, \frac{p+q-d_r-d_{r+1}}{2}-2, \ldots \frac{q-p}{2}) \\
  \preceq & (-\frac{p+q-d_r-d_{r-1}}{2}+\frac{d_{r-1}+3}{2}-\epsilon, -\frac{p+q-d_r-d_{r-1}}{2}+\frac{d_{r-1}+3}{2}-1, \ldots, \frac{-q+p}{2}+\frac{d_{r-1}+1}{2})
  \end{split}
  \end{equation}
  since $\frac{d_{r-2}}{2} \leq \frac{d_{r-1}+1}{2}- \epsilon$. Also notice that the minimal $K$-type of $ \pi_0$ is of the form $(\xi \oplus {\mathbf 0}, \eta \oplus {\mathbf 0}, +)$. By Cor. \ref{typeb1} and Remark \ref{typeb2}, the Vogan subquotient of
  $$\Ind_{SO(p-\frac{d_r+d_{r-1}}{2}, q-\frac{d_r+d_{r-1}}{2}) GL(1)^{\frac{d_r+d_{r-1}}{2}} N}^{SO(p,q)} \pi_0 \otimes ({\rm triv} \otimes \mathbb C_{ v_{(d_{r-1}, d_r)}})$$
  is unitary with a minimal $K$-type $(\xi \oplus {\mathbf 0}, \eta \oplus {\mathbf 0}, +)$. By double induction formula,
  the Langlands-Vogan parameter $(SO(p-k, q-k)GL(1)^k N, \sigma \otimes {\rm triv}, v_{\mathbf D})$ is unitary. $\Box$

\subsection{$p+q$ Odd: Type $B$ Groups}
In the case $G=SO(p,q)$ with $p+q$ odd, ${}^L G_0=Sp_{p+q-1}(\mathbb C)$. Most of the discussion here will be similar to $p+q$ even case, but not the same. We shall be brief and make the differences clear. We consider only those $\psi(w)$ a tempered parameter for $SO(p-k, q-k)$ and $\psi(SL(2, \mb C)) \subseteq Sp_{2k}(\mathbb C)$. The equivalence classes of homomorphisms from $SL(2, \mathbb C)$ to $Sp_{2k}(\mathbb C)$ are in one-to-one correspondence with nilpotent adjoint orbits of $Sp_{2k}(\mathbb C)$, which is in one-to-one correspondence to the symplectic Young diagrams of size $2k$. A Young diagram is called symplectic if
every odd part occurs with even multiplicity. \\
\\
Let $\mathbf D$ be a symplectic Young diagram of $2k$. It can be described as a partition of 
$2k$:
$$\sum_{i=1}^r d_i=2k, \qquad (d_1 \leq d_2 \ldots d_{r-1} \leq d_r).$$
Let $s_{\mathbf D}$ be the semisimple element corresponding to ${\rm diag}(\frac{1}{2}, -\frac{1}{2})$. Use $v_{\mathbf D} \in \mathbb R^{k}$ to parametrize the conjugacy class of $s_{\mathbf D}$. Then
$v_{\mathbf D}$ will be unique up to sign changes and permutations. $v_{\mathbf D}$ can be constructed as a direct sum in the following way (Pg 78 \cite{cm}).
\begin{enumerate}
\item If $d_i=d_{i+1}=d$, then $v_{(d, d)}=(\frac{d-1}{2}, \frac{d-3}{2}, \ldots, \frac{1}{2}, -\frac{1}{2}, \ldots, \frac{1-d}{2})$. We do this for every pair of integers. Then we are only left with even $d_i$'s.
\item If $d_i$ is even, let $v_{(d_i)}=(\frac{d_i-1}{2}, \frac{d_i-3}{2}, \ldots, \frac{1}{2})$. We do this for every even $d_i$ that is not in an identical pair.
\end{enumerate} 

Now we can restate Theorem B for $p+q$ odd.
 \begin{cor}\label{typed1}
 Suppose that $p+q \leq 2m+1 \leq p+q+d$, $p \leq q$ and $p+q$ odd. Let $\pi$ be an irreducible unitary representation of $SO(p,q)$ such that
its every $K$ finite matrix coefficient $f(g)$ satisfies the condition that
$$ |f(g)| \leq C_f \exp (\overbrace{m+2-p-q- \epsilon, m+3-p-q, \ldots m+1-q}^{p})(| H^+(g)|).$$
for some $\epsilon >0$. 
Suppose that the minimal $S(O(p) O(q))$-types of $\pi$ contains
a $(\xi, \eta, +)$. Then the Vogan subquotient of 
$\Ind_{SO(p,q) GL(1)^{d} N}^{SO(p+d, q+d)} \pi \otimes ({\rm triv} \otimes \mathbb C_{ v_{(2m+1-p-q)} \oplus v_{(p+q+2d-2m-1)}})$
is unitary.
 \end{cor}
 \begin{re}\label{typed2} Now let $2m+1-p-q=s$ and $p+q+2d-2m-1=t$. Then $s,t$ are even. The growth condition becomes
 $$ |f(g)| \leq C_f \exp (\overbrace{-\frac{p+q}{2}+\frac{s+3}{2}- \epsilon, -\frac{p+q}{2}+\frac{s+5}{2}, \ldots \frac{-q+p}{2}+\frac{s+1}{2}}^{p})(| H^+(g)|).$$
 The assertion is that the Vogan subquotient of $\Ind_{SO(p,q) GL(1)^{d} N}^{SO(p+d, q+d)} \pi \otimes ({\rm triv} \otimes \mathbb C_{ v_{(s)} \oplus v_{(t)}})$ is unitary.
 \end{re}

\begin{Theorem}[Theorem D] Let $p+q$ be odd. Let $\mathbf D$ be an symplectic Young diagram of size $2k$. Let $\sigma$ be an irreducible tempered representation of $SO(p-k, q-k)$ such that $\sigma$ has a minimal $K$-type of the form $(\xi, \eta, +)$. Then the Langlands-Vogan parameter $(SO(p-k, q-k)GL(1)^k N, \sigma \otimes {\rm triv}, v_{\mathbf D})$ is unitary,  with a minimal $K$-type
 $(\xi \oplus {\mathbf 0}, \eta \oplus {\mathbf 0}, +)$.
 \end{Theorem}
 Proof: Clearly,  a minimal $K$-type of  $(SO(p-k, q-k)GL(1)^k N, \sigma \otimes {\rm triv}, v_{\mathbf D})$  is of the form 
 $(\xi \oplus {\mathbf 0}, \eta \oplus {\mathbf 0}, +)$. Let $\mathbf D$ be defined by the partition
 $$d_1 \leq d_2 \leq d_3 \ldots \leq d_r \qquad (\sum_{i=1}^r d_i=2k).$$ 
 We apply induction on $r$. If $r=0$, we have the tempered representation $\sigma$ which is unitarizable.\\
 \\
 If there is a pair $d_i=d_{i+1}=d$, let $\mathbf D_0$ be the Young diagram obtained by deleting $d_i, d_{i+1}$. Then by induction hypothesis, the Langlands-Vogan parameter $(SO(p-k, q-k)GL(k-d)N, \sigma \otimes {\rm triv}, v_{\mathbf D_0})$ is unitary. Denote this representation by $\pi_0$. Then $\Ind_{SO(p-d, q-d) GL(d) N}^{SO(p,q)} \pi_0 \otimes {\rm triv}$ is unitary. It is easy to see that the Langlands-Vogan parameter $(SO(p-k,q-k) GL(1)^k N, \sigma \otimes {\rm triv}, v_{\mathbf D})$ is a subrepresentation of 
 $\Ind_{SO(p-d, q-d) GL(d) N}^{SO(p,q)} \pi_0 \otimes {\rm triv}$, hence is unitary. \\
 \\
 If there are no repeats in $d_i$, then $d_1, d_2, \ldots, d_r$ must be all even. Now pair all the $d_i$'s in descending order, namely, $(d_r, d_{r-1})$, $(d_{r-2}, d_{r-3})$ and so on. Add $d_0=0$ if necessary.
 Let $\mathbf D_0$ be the Young diagram obtained from $\mathbf D$ by deleting $d_r$ and $d_{r-1}$.  
 By induction hypothesis, the Langlands-Vogan parameter $(SO(p-k,q-k) GL(1)^{k-\frac{d_r+d_{r-1}}{2}} N, \sigma \otimes {\rm triv}, v_{\mathbf D_0})$ is unitary. Denote this unitary representation of $SO(p-\frac{d_r+d_{r-1}}{2}, q-\frac{d_r+d_{r-1}}{2})$ by $\pi_0$. Then by  Remark \ref{langlandsmatrixco}, $K$-finite matrix coefficients of $\pi_0$ are bounded by multiples of $\Xi_{v_{\mathbf D_0} \oplus {\bf 0_{p-k}}}(g)$. The leading exponents of $\Xi_{v_{\mathbf D_0} \oplus {\bf 0_{p-k}}}(g)$ is bounded by
 \begin{equation}
 \begin{split}
  & (\frac{d_{r-2}}{2}, \frac{d_{r-2}}{2}, \ldots, \frac{d_{r-2}}{2}) -(\frac{p+q-d_r-d_{r-1}}{2}-1, \frac{p+q-d_r-d_{r+1}}{2}-2, \ldots \frac{q-p}{2}) \\
  \preceq & (-\frac{p+q-d_r-d_{r-1}}{2}+\frac{d_{r-1}+3}{2}-\epsilon, -\frac{p+q-d_r-d_{r-1}}{2}+\frac{d_{r-1}+3}{2}-1, \ldots, \frac{-q+p}{2}+\frac{d_{r-1}+1}{2}) \\
  \end{split}
  \end{equation}
  since $\frac{d_{r-2}}{2}+1 \leq \frac{d_{r-1}+3}{2}- \epsilon$. Also notice that the minimal $K$-type of $\pi_0$ is of the form $(\xi \oplus {\mathbf 0}, \eta \oplus {\mathbf 0}, +)$. By Cor. \ref{typed1} and Remark \ref{typed2}, the Vogan subquotient of  
  $$\Ind_{SO(p-\frac{d_r+d_{r-1}}{2}, q-\frac{d_r+d_{r-1}}{2}) GL(1)^{\frac{d_r+d_{r-1}}{2}} N}^{SO(p,q)} \pi_0 \otimes ({\rm triv} \otimes \mathbb C_{ v_{(d_{r-1})} \oplus v_{(d_r)}})$$
  is unitary. By double induction formula,
  the Langlands-Vogan parameter $$(SO(p-k, q-k)GL(1)^k N, \sigma \otimes {\rm triv}, v_{\mathbf D})$$ is unitary.  $\Box$

\end{document}